\documentclass[12pt, leqno]{article}
\usepackage[final]{epsfig}
\usepackage{graphics}
\usepackage{amsmath}
\usepackage{amsfonts}
\usepackage{latexsym}
\usepackage{amssymb}
\usepackage{graphicx}
\usepackage{url}
\usepackage{epstopdf}
\usepackage{xcolor}
\def\Lie{\mathop{\rm Lie}}
\def\Stab{\mathop{\rm Stab}}
\def\Ind{\mathop{\rm Ind}}
\def\tr{\mathop{\rm tr}}
\input pictexwd.tex

\font\small=cmr7
\font\huge=cmr24

\newtheorem{theorem}{Theorem}[section]
\newtheorem{lemma}[theorem]{Lemma}
\newtheorem{proposition}[theorem]{Proposition}
\newtheorem{corollary}[theorem]{Corollary}

\begin{document}

\title{Representations of the group of two-diagonal triangular matrices.}

\author{Dmitry Fuchs\footnote{Department of Mathematics, University of California, Davis,
 fuchs@math.ucdavis.edu} \and Alexandre Kirillov\footnote{Department of Mathematics, University of Pennsylvania, kirillov@math.upenn.edu}
}

\maketitle 

\hskip1.6in{\it To Leonid Makar-Limanov on the ocasion of his 80-th birthday}

\begin{abstract} {\rm Let $G$ be a Lie group, ${\mathfrak g} = \Lie(G)$ – its Lie algebra, ${\mathfrak g}^\ast$ – the dual vector space and $\widehat G$ – the set of equivalence classes of unitary irreducible representations of $G$. The orbit method [1] establishes a correspondence between points of $\widehat G$ and $G$-orbits in ${\mathfrak g}^\ast$. For many Lie groups it gives the answers to all major problems of representation theory in terms of coadjoint orbits. Formally, the notions and statements of the orbit method make sense when $G$ is infinite-dimensional Lie group, or an algebraic group over a topological field or ring $K$, whose additive group is self dual (e.g., $p$-adic or finite).

In this paper, we introduce the big family of finite groups $G_n$, for which the orbit method 
works perfectly well. Namely, let $N_n({\mathbb K})$ be the algebraic group of upper unitriangular 
$(n+1)\times(n+1)$ matrices with entries from ${\mathbb K}$, and ${\mathbb F}_q$ be the finite field with $q$ elements. We define $G_n$ as the quotient of of the group $N_{n+1}({\mathbb F}_q)$ over its second commutator subgroup.}\end{abstract}

\section{Introduction}
This paper is a part of a bigger program of the application of the orbit method in the representation theory. The main ingredient of the orbit method is the notion of a {\bf coadjoint orbit}. The method works not only for ordinary Lie groups, but also for infinite-dimensional Lie groups and for algebraic groups over a topological ring $\mathbb K$, which satisfies the following condition. {\it Every additive character of $\mathbb K$ has the form $\chi_\lambda(x)={\bf e}(\lambda x)$, where $\lambda\in{\mathbb K}$ and $\bf e$ is a fixed non-trivial character.}\footnote{In other words, the additive group ${\mathbb K}^+$ is Pontryagin self-dual.} This condition is satisfied for real, complex, quaternion, $p$-adic and finite fields and also for the adel ring.

It is known (see. e.g. \cite{kir}) that the method can be formulated as the collection of simple rules, which gives the transparent answers to all main questions. This ``User's guide" can be understood literally when the group in question is a connected and simply connected nilpotent Lie group. For many other Lie or algebraic groups it works after some algebraic and/or topological corrections. There are also groups, for which the right application of the orbit method is still unknown.

One time it was a hope that some modification of the orbit method will work for the groups $N_n({\mathbb F}_q)$ of upper unitriangular $(n+1)\times(n+1)$ matrices with elements from a finite field ${\mathbb F}_q$ (considered as algebraic groups). A simple argument shows that the number of coadjoint orbits of the group $N_n({\mathbb F}_q)$ is equal to the number of conjugacy classes, hence to the number of equivalence classes of irreducible representations of this group. However, the explicit construction of irreducible representations and their characters in terms of coadjoint orbits is known only for small $n$. 

In this paper, we consider the family of groups related to, but more simple than the full triangular group $N_n({\mathbb K})$, for which all questions of representation theory have explicit answers in terms of coadjoint orbits. Namely, we consider the quotient group of $N_n({\mathbb K})$ by its second commutant $[[N_n({\mathbb K}),N_n({\mathbb K})],N_n({\mathbb K})]$. We call it ``two-diagonal group" and denote by $TD_n({\mathbb K})$, or simply $G_n$. It is a $(2n-1)$-dimensional affine algebraic group over  the field $\mathbb K$. As such the group $G_n$ has the Lie algebra ${\mathfrak g}_n$ with the dual space ${\mathfrak g}_n^\ast$. These spaces are the spaces of respectively, {\it adjoint} and {\it coadjoint} representations of $G_n$. Our first goal is to describe explicitly orbits of these two representations. 

At the beginning, we impose no restrictions on the ground field $\mathbb K$, but out main results (Sections 2.2, 3.2, 4, 5) concern the case $\mathbb{K=F}_q$. 

Both ${\mathfrak g}_n$ and ${\mathfrak g}^\ast_n$ can be realized as subquotients of the full matrix space Mat$_n(\mathbb K)$. We choose the coordinates $\{a_i\}_{1\le i\le n},\{b_j\}_{1\le j\le n-1}$ for $A\in{\mathfrak g}_n$ and $\{x_i\}_{1\le i\le n},\{y_j\}_{1\le j\le n-1}$ for $F\in{\mathfrak g}^\ast_n$, so that
$$A=\left[\begin{array} {cccccc} 0&a_1&b_1&\dots&0&0\\ 0&0&a_2&\dots&0&0\\ \dots&\dots&\dots&\dots&\dots&\dots\\ 0&0&0&\dots&a_{n-1}&b_{n-1}\\ 0&0&0&\dots&0&a_n\\ 0&0&0&\dots&0&0\end{array}\right],F=\left[\begin{array} {cccccc} 0&0&\dots&0&0&0\\ x_1&0&\dots&0&0&0\\ y_1&x_2&\dots&0&0&0\\ \dots&\dots&\dots&\dots&\dots&\dots\\ 0&0&\dots&x_{n-1}&0&0\\ 0&0&\dots&y_{n-1}&x_n&0\end{array}\right].$$A general element of the group $G_n$ has the form $$\left[\begin{array} {cccccc} 1&\alpha_1&\beta_1&\dots&0&0\\ 0&1&\alpha_2&\dots&0&0\\ \dots&\dots&\dots&\dots&\dots&\dots\\ 0&0&0&\dots&\alpha_{n-1}&\beta_{n-1}\\ 0&0&0&\dots&1&\alpha_n\\ 0&0&0&\dots&0&1\end{array}\right]$$ This justifies the name ``two-diagonal group"). We denote this element by $g(\overline\alpha;\overline\beta)$ where $\overline\alpha=(\alpha_1,\dots,\alpha_n)$ and $\overline \beta=(\beta_1,\dots,\beta_{n-1})$. The  group laws are $$\begin{array} {c} g(\alpha_1,\dots.,\alpha_n;\beta_1,\dots,\beta_{n-1})g(\alpha'_1,\dots.,\alpha'_n;\beta'_1,\dots,\beta'_{n-1})=\hskip1.9in\\ \hskip.8in g(\alpha_1+\alpha'_1,\dots,\alpha_n+\alpha'_n;\beta_1+\beta'_1+\alpha_1\alpha'_2,\dots,\beta_{n-1}+\beta'_{n-1}+\alpha_{n-1}\alpha'_n);\\ g(\alpha_1,\dots.,\alpha_n;\beta_1,\dots,\beta_{n-1})^{-1}=g(-\alpha_1,\dots,-\alpha_n;-\beta_1+\alpha_1\alpha_2,\dots,-\beta_{n-1}+\alpha_{n-1}\alpha_n)\end{array}$$

For adjoint and coadjoint representations, the action of $g(\overline\alpha,\overline\beta)\in G_n$ does not depend on $\overline\beta$. The diagonals of $a$'s and $y$'s are invariant with respect to the $G_n$-actions, and the action on the diagonals of $b$'s and $x$'s is described by the formulas: 

$$\left[\hskip-.07in\begin{array} {c} b'_1\\b'_2\\ \dots\\ b'_{n-1}\end{array}\hskip-.07in\right]=\left[\hskip-.07in\begin{array} {c} b_1\\b_2\\ \dots\\ b_{n-1}\end{array}\hskip-.07in\right]+[\alpha_1\, \alpha_2\, \dots\, \alpha_n]\cdot S,$$\vskip-.2in $$[x'_1\, x'_2\, \dots\, x'_n]=[x_1\, x_2\, \dots\, x_n]+\left[\hskip-.07in\begin{array} {c} \alpha_1\\ \alpha_2\\ \dots\\ \alpha_n\end{array}\hskip-.07in\right]\cdot T,$$where the $(n-1)\times n$ matrix $S$ and the $n\times n$ matrix $T$ are
$$S=\left[\hskip-.05in\begin{array} {rrrrc} a_2&0&\dots&0&0\\ -a_1&a_3&\dots&0&0\\ \dots&\dots&\dots&\dots&\dots\\ 0&0&\dots&a_{n-1}&0\\ 0&0&\dots&-a_{n-2}&a_n\\ 0&0&\dots&0&-a_{n-1}\end{array}\hskip-.05in\right],\ T=\left[\hskip-.05in\begin{array} {rrrrccc} 0&y_1&0&\dots&0&0&0\\ -y_1&0&y_2&\dots&0&0&0\\ 0&-y_2&0&\dots&0&0&0\\ \dots&\dots&\dots&\dots&\dots&\dots&\dots\\ 0&0&0&\dots&0&y_{n-2}&0\\ 0&0&0&\dots&-y_{n-2}&0&y_{n-1}\\ 0&0&0&\dots&0&-y_{n-1}&0\end{array}\right].$$ It will be more convenient to us to write the above formulas for the $G_n$-action in a more traditional form:

$$\begin{array} {ccccc} b'_1=b_1&&+\alpha_1a_2&&-\alpha_2a_1\\ b'_2=b_2&&+\alpha_2a_3&&-\alpha_3a_2\\ b'_3=b_3&&+\alpha_3a_4&&-\alpha_4a_3\end{array} \hskip.6in \begin{array} {ccccc} x'_1=x_1&&&&-\alpha_2y_1\\ x'_2=x_2&&+\alpha_1y_1&&-\alpha_3y_2\\ x'_3=x_3&&+\alpha_2y_2&&-\alpha_4y_3\end{array}$$\vskip-.14in

$$\dots\dots\dots\dots\dots\dots\dots\dots\dots\hskip.8in\dots\dots\dots\dots\dots\dots\dots\dots\dots\dots\eqno(1)$$\vskip-.3in

$$\begin{array} {l} b'_{n-2}=b_{n-2}+\alpha_{n-2}a_{n-1}-\alpha_{n-1}a_{n-2}\\  b'_{n-1}=b_{n-1}+\alpha_{n-1}a_n-\alpha_na_{n-1}\end{array}\hskip.3in\begin{array} {l} x'_{n-1}=x_{n-1}+\alpha_{n-2}y_{n-2}-\alpha_ny_{n-1}\\ x'_n=x_n+\alpha_{n-1}y_{n-1}\end{array}$$

We will refer to orbits of the coadjoint representation as just to ``orbits," and to orbits of the adjoint representation as to ``classes." 

\section{Orbits.}

\subsection{Description.} 

For given $y_1,\dots,y_{n-1}; x_1,\dots,x_n; x'_1,\dots,x'_n$ we consider the right system (1) as the system of $n$ equations with $n$ unknowns $\alpha_1,\dots\alpha_n$. If this system has a solution, then the matrices with $x_1,\dots,x_n$ and $x'_1,\dots,x'_n$  belong to the same orbit.

First consider the case, when $y_1,\dots,y_{n-1}$ are all different from zero. Our system splits into two systems with unknowns $\alpha_1,\alpha_3,\alpha_5\dots$ and unknowns $\alpha_2,\alpha_4,\alpha_6,\dots$, and this splitting looks differently in the cases of even and odd $n$.\smallskip

{\bf If n is even}, then the two systems are

$$\begin{array} {ccccc} x'_2=x_2&&+\alpha_1y_1&&-\alpha_3y_2\\ x'_4=x_4&&+\alpha_3y_3&&-\alpha_5y_4 \end{array} \hskip.6in \begin{array} {ccccc} x'_1=x_1&&&&-\alpha_2y_1\\ x'_3=x_3&&+\alpha_2y_2&&-\alpha_4y_3 \end{array}$$\vskip-.14in

$$\dots\dots\dots\dots\dots\dots\dots\dots\dots\dots\quad{\rm and}\quad\dots\dots\dots\dots\dots\dots\dots\dots\dots\dots$$\vskip-.18in

$$\begin{array} {l} x'_{n-2}=x_{n-2}+\alpha_{n-3}y_{n-3}-\alpha_{n-1}y_{n-2}\\  x'_n=x_n+\alpha_{n-1}y_{n-1} \end{array}\hskip.3in\begin{array} {l} x'_{n-3}=x_{n-3}+\alpha_{n-4}y_{n-4}-\alpha_{n-2}y_{n-3}\\ x'_{n-1}=x_{n-1}+\alpha_{n-2}y_{n-2}-\alpha_ny_{n-1}\end{array}$$

Both systems have unique solutions. We solve the first system ``from the bottom to the top": we find $\alpha_{n-1}$ from the last equation, then we find $\alpha_{n-3}$ from the equation before the last, and so on, up to $\alpha_3$ from the second equation and $\alpha_1$ from the first equation. The second system can be solved ``from the top to the bottom": we find $\alpha_2$ from the first equation, then find $\alpha_4$ from the second equation, and so on. Thus, in this case there is only one orbit (corresponding to the fixed set $y_1,\dots,y_{n-1}$ of non-zero elements of $\mathbb K$). This orbit is an $n$-dimensional vector space over ${\mathbb K}$ with coordinates $x_1,\dots,x_n.$ \smallskip

{\bf If n is odd}, then the two systems are

$$\hskip.6in \begin{array} {ccccc} x'_1=x_1&&&&-\alpha_2y_1\\ x'_3=x_3&&+\alpha_2y_2&&-\alpha_4y_3 \end{array}\hskip.9in\begin{array} {ccccc} x'_2=x_2&&+\alpha_1y_1&&-\alpha_3y_2\\ x'_4=x_4&&+\alpha_3y_3&&-\alpha_5y_4 \end{array} $$\vskip-.14in

$$\hskip.2in\dots\dots\dots\dots\dots\dots\dots\dots\dots\dots\quad{\rm and}\quad\dots\dots\dots\dots\dots\dots\dots\dots\dots\dots$$\vskip-.18in

$$\hskip.2in\begin{array} {l} x'_{n-2}=x_{n-2}+\alpha_{n-3}y_{n-3}-\alpha_{n-1}y_{n-2}\\  x'_n=x_n+\alpha_{n-1}y_{n-1} \end{array}\hskip.3in\begin{array} {l} x'_{n-3}=x_{n-3}+\alpha_{n-4}y_{n-4}-\alpha_{n-2}y_{n-3}\\ x'_{n-1}=x_{n-1}+\alpha_{n-2}y_{n-2}-\alpha_ny_{n-1}\end{array}$$

For the second system, the number of equations is one less than the number of unknowns. It has more than one solution: we can choose arbitrary $\alpha_n$ and then find  $\alpha_{n-2},\alpha_{n-4},\dots,\alpha_3,\alpha_1$ from our equations from the bottom to the top. However, for the first system, the number of equations is one more than the number of unknown. We can find $\alpha_2,\alpha_4,\dots,\alpha_{n-3},\alpha_{n-1}$ using all the equations except the last one, $x'_n=x_n+\alpha_{n-1}y_{n-1}$, and the value of $\alpha_{n-1}$, which we have already found, may unfit this last equation. A simple computation shows that this equation is satisfied is the equality $$\begin{array} {l}\hskip.14in  x_1y_2y_4\dots y_{n-1}+y_1x_3y_4y_6\dots y_{n-1}+y_1y_3x_5y_6y_8\dots y_{n-1}\\ \hskip1.8in+\dots+y_1y_3\dots y_{n-4}x_{n-2}y_{n-1}+y_1y_3\dots y_{n-2}x_n=\\ =x'_1y_2y_4\dots y_{n-1}+y_1x'_3y_4y_6\dots y_{n-1}+y_1y_3x'_5y_6y_8\dots y_{n-1}\\ \hskip1.8in+\dots+y_1y_3\dots y_{n-4}x'_{n-2}y_{n-1}+y_1y_3\dots y_{n-2}x'_n.\end{array}\eqno(2)$$holds. In other words, the expression in the left hand side of (2) is an invariant of the action of $G_n$, and, {\bf in the case of odd $\bf n$}, the orbit is characterized by the set $y_1,\dots,y_{n-2}$ of non-zero elements of $\mathbb K$ and the value $I\in{\mathbb K}$ of our invariant (which may be 0).The dimension of the orbit in this case is $n-1$ with coordinates $x_1,\dots,x_{n-1}.$\smallskip

Now let us turn to the general case: suppose that $y_{i_1}=\dots=y_{i_m}=0\ (1\le i_1<\dots<i_m\le n-1)$, while all the remaining $n-1-m$ $y$'s are different from zero. Then the right system (1) splits into $m+1$ independent systems, $\bf S_0,S_1,\dots,S_m$:

\centerline{
\beginpicture
\setcoordinatesystem units <1in,1in> point at 0 0
\put{System $\bf S_0$} at -2.2 -.35
\put{System $\bf S_1$} at -2.2 -1.35
\setplotsymbol({\small.})
\put{$x'_1=x_1$} at -1 0
\put{$x'_2=x_2$} at -1 -.2
\put{$\dots\dots\dots\dots\dots\dots\dots\dots\dots\dots\dots\dots\dots\dots$} at -.2 -.35
\put{$x'_{i_1-1}=x_{i_1-1}$} at -1 -.5
\put{$x'_{i_1}=x_{i_1}$} at -1 -.7
\putrule from -2.45 -.85 to 1.1 -.85
\put{$x'_{i_1+1}=x_{i_1+1}$} at -1 -1
\put{$x'_{i_1+2}=x_{i_1+2}$} at -1 -1.2
\put{$\dots\dots\dots\dots\dots\dots\dots\dots\dots\dots\dots\dots\dots\dots$} at -.2 -1.35
\put{$x'_{i_2-1}=x_{i_2-1}$} at -1 -1.5
\put{$x'_{i_2}=x_{i_2}$} at -1 -1.7
\putrule from -2.45 -1.85 to 1.1 -1.85
\put{$\dots\dots\dots\dots\dots\dots\dots\dots\dots\dots\dots\dots\dots\dots\dots\dots\dots
                          \dots\dots\dots\dots\dots$} at -.4 -2
\put{$+\alpha_1y_1$} at  -.1 -.2
\put{$+\alpha_{i_1-2}y_{i_1-2}$} at -.1 -.5
\put{$+\alpha_{i_1-1}y_{i_1-1}$} at -.1 -.7
\put{$+\alpha_{i_1+1}y_{i_1+1}$} at -.1 -1.2
\put{$+\alpha_{i_2-2}y_{i_2-2}$} at -.1 -1.5
\put{$+\alpha_{i_2-1}y_{i_2-1}$} at -.1 -1.7
\put{$-\alpha_2y_1$} at .7 0
\put{$-\alpha_3y_2$} at .7 -.2
\put{$-\alpha_{i_1}y_{i_1-1}$} at .7 -.5
\put{$-\alpha_{i_1+2}y_{i_1+1}$} at .7 -1
\put{$-\alpha_{i_1+3}y_{i_1+2}$} at .7 -1.2
\put{$-\alpha_{i_2}y_{i_2-1}$} at .7 -1.5
\put{$(y_{i_1}=0)$} at 1.45 -.85
\put{$(y_{i_2}=0)$} at 1.45 -1.85
\endpicture
}\vskip.2in

\centerline{
\beginpicture
\setcoordinatesystem units <1in,1in> point at 0 0
\put{$\dots\dots\dots\dots\dots\dots\dots\dots\dots\dots\dots\dots\dots\dots\dots\dots\dots
                          \dots\dots\dots\dots\dots$} at -.4 -2
\putrule from -2.45 -2.15 to 1.1 -2.15
\put{System $\bf S_m$} at -2.2 -2.65
\put{$x'_{i_m+1}=x_{i_m+1}$} at -1 -2.3
\put{$x'_{i_m+2}=x_{i_m+2}$} at -1 -2.5
\put{$\dots\dots\dots\dots\dots\dots\dots\dots\dots\dots\dots\dots\dots\dots$} at -.2 -2.65
\put{$x'_{n-1}=x_{n-1}$} at -1 -2.8
\put{$x'_n=x_n$} at -1 -3
\put{$+\alpha_{i_m+1}y_{i_m+1}$} at -.04 -2.5
\put{$+\alpha_{n-2}y_{n-2}$} at -.04 -2.8
\put{$+\alpha_{n-1}y_{n-1}$} at -.04 -3
\put{$-\alpha_{i_1+2}y_{i_1+1}$} at .8 -2.3
\put{$-\alpha_{i_1+3}y_{i_1+2}$} at .8 -2.5
\put{$-\alpha_ny_{n-1}$} at .8 -2.8
\put{$(y_{i_m}=0)$} at 1.45 -2.15
\endpicture
}\vskip.2in

Let us introduce additional notations: $i_0=0,i_{m+1}=n$, and $j_r=i_{r+1}-i_r$ for $r=0,1,\dots,k$. Thus, $j_0+j_1+\dots+j_m=n$.

The systems $\bf S_0,S_1,\dots,S_m$ have disjoint sets of unknowns: the $j_r$ unknowns in the system $\bf S_r$ are $\alpha_{i_r+1},\dots,\alpha_{i_{r+1}}$ (up to the numerations of $\alpha$'s, $x$'s $x'$'s, and $y$'s) are reduced copies of the right system (1), and we can apply to it our findings for that system. Thus, if $j_r$ is even, then the system $\bf S_r$ has a unique solution $\alpha_{i_r+1},\dots,\alpha_{i_{r+1}}$, and if $j_r$ is odd then the system $\bf S_r$ is consistent (although a solution is not unique) if and only if $I_r=I'_r$, where$$\begin{array} {l} I_r=x_{i_r+1}y_{i_r+2}y_{i_r+4}\dots y_{i_{r+1}-1}+y_{i_r+1}x_{i_r+3}y_{i_r+4}y_{i_r+6}\dots y_{i_{r+1}-1}+\dots+\\ \hskip1in+y_{i_r+1}y_{i_r+3}\dots y_{i_{r+1}-5}x_{i_{r+1}-3}y_{i_{r+1}-2}+y_{i_r+1}y_{i_r+3}\dots y_{i_{r+1}-3}x_{i_{r+1}-1},\end{array}$$and $I'_r$ is defined by the same formula with all $x$'s replaced by $x'$'s.\smallskip

We arrive at the following description of all orbits of $G_n$. 

First, we need to fix an ordered  partition$$n=j_0+j_1+\dots+j_m$$(where $j_0,j_1,\dots,j_m$ are positive integers). Then we put $$\begin{array} {l} i_1=j_0,\\ i_2=j_0+j_1,\\ i_3=j_0+j_1+j_2,\\ \dots\dots\dots\dots\dots\\ i_m=j_0+j_1+\dots+j_{m-1}.\end{array}$$

Second, we choose $y_1,\dots,y_{n-1}\in\mathbb K$ such that $y_{i_1}=\dots=y_{i_m}=0$ and all the other $y$'s are different from zero. 

Third, we choose a $v_r\in\mathbb K$ for each $r$ from $0,1,\dots,m$ such that $j_r$ is odd.

These data determine an orbit. This orbit consists of the matrices with $y$'s fixed above and arbitrary $x_1,\dots,x_n\in\mathbb K$ satisfying the condition $I_r=v_r$ for all $r$ such that $j_r$ is odd. If $\nu$ is the number of such $r$ then the dimension of this orbit is $n-\nu$. Notice that $\nu$ and $n$ have the same parity, $n-\nu=2k$ for an integer $k,\, 0\leq k\leq \left[\displaystyle\frac n2\right]$. Hence the dimension of the orbit is $2k$, in particular, {\bf the dimensions of all orbits are even}.

If ${\mathbb K=\mathbb F}_q$, then the number of orbits corresponding to an ordered partition $n=j_0+j_1+\dots+j_m$ is $(q-1)^{n-m-1}q^\nu$. 

Notice that if $\mu$ is the number of {\it even} terms of our partition, then $\mu\leq k$. Indeed, each even term of the partition is at least 2, and each odd term is at least 2. Hence $n\geq 2\mu+\nu=2\mu+(n-2k)=n-2(k-\mu)\Rightarrow k-\mu\geq0$.

\subsection{The number of orbits of given dimension.}

 We suppose again that ${\mathbb K=F}_q$. The following holds:\smallskip

\begin{theorem}The  number of orbits of $G_n$ of dimension $2k$ is $$q^{n-k-1}(q-1)^k\left({n-k-1\choose k}q+{n-k-1\choose k-1}\right).$$\end{theorem}

The proof is based on the following combinatorial\smallskip

\begin{lemma} The number of (ordered) partition of $n\geq2\mu+\nu$ into the sum of $\mu$ even and $\nu$ odd positive integers is $${\mu+\nu\choose \mu}\cdot{(n+\nu)/2-1\choose \mu+\nu-1}.$$\end{lemma}

\noindent{\bf Proof of Lemma 2.2.} To specify a partition of $n$ into $\mu$ even and $\nu$ odd summands, we first need to choose $\mu$ positions for the even summands, which can be done in $\displaystyle{\mu+\nu\choose \mu}$ ways. Then we transform all our partitions in the following way: we add 1 to each odd summand and then divide all the summands by 2. We get an ordered partition of the number $\displaystyle\frac{n+\nu}2$ into $\mu+\nu$ summands with no condition on the parity of the summands. It is well known that the number of such partitions is $\displaystyle{(n+\nu)/2-1\choose \mu+\nu-1}$, which implies our statement.\smallskip

\noindent{\bf Proof of Theorem 2.1.}  According to Section 2.1, $2k$-dimensional orbits of $G_n$ correspond to (ordered) partitions of $n$ with the number of odd terms equal to $\nu=n-2k$. The number of orbits corresponding to partitions with $\mu$ even and $\nu$ odd parts is$$(q-1)^{n-(\mu+\nu)}q^\nu=(q-1)^{2k-\mu}q^{n-2k}.$$By Lemma, the number of such partitions is $$\begin{array} {l} \displaystyle{{\mu+\nu\choose \mu}{(n+\nu)/2-1\choose \mu+\nu-1}={n-2k+\mu\choose\mu}{n-k-1\choose k-\mu}}=\\ \\ \displaystyle{\frac{(n-k-1)(n-k-2)\dots(n-2k+\mu)}{(k-\mu)!}\cdot\frac{(n-2k+\mu)\dots(n-2k+1)} {\mu!}}=\\ \\ \hskip1.6in\displaystyle\frac{(n-k-1)(n-k-2)\dots(n-2k+1)}{(k-\mu)!\mu!}\cdot(n-2k+\mu).\end{array}$$

Thus the total number of orbits of dimension $2k$ is $$\begin{array} {l} \displaystyle{\sum_{b=0}^k\frac{(n-k-1)(n-k-2)\dots(n-2k+1)}{(k-\mu)!\mu!}}(q-1)^{2k-\mu}q^{n-2k}\cdot(n-2k+\mu)\\ \displaystyle{=\frac{(n-k-1)\dots(n-2k+1)}{k!}q^{n-2k}(q-1)^k\sum_{\mu=0}^k(n-2k+\mu){k\choose \mu}(q-1)^{k-\mu}}\end{array}$$The sum in the last line is$$\begin{array} {c} \displaystyle{(n-k)\sum_{\mu=0}^k{k\choose \mu}(q-1)^{k-\mu}-\sum_{b=0}^k(k-\mu){k\choose \mu}(q-1)^{k-\mu}=(n-2k)q^k-(q-1)(q^k)'}\\ \\ =(n-k)q^k-k(q-1)q^{k-1}=(n-2k)q^k+kq^{k-1}=q^{k-1}((n-2k)q+k).\end{array}$$

We plug this expression for the sum into the last formula and see that the total number of orbits of dimension $2k$ is$$\begin{array} {c} \displaystyle{q^{n-k-1}(q-1)^k\left(\frac{(n-k-1)\dots(n-2k-1)}{k!}q+\frac{(n-k-1)\dots(n-2k)}{(k-1)!}\right)}\\ \\ \displaystyle{=q^{n-k-1}(q-1)^k\left({n-k-1\choose k}q+{n-k-1\choose k-1}\right).}\end{array}$$

\subsection{More on ordered partitions.} 

Since ordered partitions of natural numbers played an essential role in the previous section, we present here a survey of their properties. Technically, the results of these section will be used only in Section 5 below, so the reader may postpone its reading until that section.

\subsubsection{Partial ordering. Intervals.} There is a partial ordering in the set of all ordered partition of a number $n$: a partition ${\mathcal P}_1$ of $n$ {\it precedes} a partition ${\mathcal P}_2$ (notations: ${\mathcal P}_1\preceq{\mathcal P}_2$ or ${\mathcal P}_2\succeq{\mathcal P}_1$), if the parts of ${\mathcal P}_2$ are obtained from the parts of ${\mathcal P}_1$ by a further subdivision. Notice that if $\nu_i$ is the number of odd terms of the partition ${\mathcal P}_i$ and ${\mathcal P}_1\preceq{\mathcal P}_2$, then $\nu_1\leq\nu_2$ (there is no  such inequality for even terms).

A more visualizable description of this partial ordering can be done in terms of dividers. Namely, an ordered partition of $n$ may be presented as a line of $n$ dots with dividers placed between some of them. For example, the diagram \vskip.15in

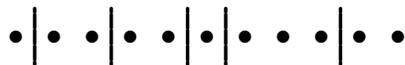
\centerline{
\beginpicture
\setcoordinatesystem units <1in,1in> point at 0 0
\setplotsymbol(.)
\put{$\bullet$} at 1 0
\put{$\bullet$} at .8 0
\put{$\bullet$} at .6 0
\put{$\bullet$} at .4 0
\put{$\bullet$} at .2 0
\put{$\bullet$} at 0 0
\put{$\bullet$} at -.2 0
\put{$\bullet$} at -.4 0
\put{$\bullet$} at -.6 0
\put{$\bullet$} at -.8 0
\put{$\bullet$} at -1 0
\plot .7 .15 .7 -.15 /
\plot .1 .15 .1 -.15 /
\plot -.1 .15 -.1 -.15 /
\plot -.5 .15 -.5 -.15 /
\plot -.9 .15 -.9 -.15 /
\endpicture
}  \vskip.25in

\noindent presents the partition $1+2+2+1+3+2$ of the number 11. In the language of such diagrams, the relation ${\mathcal P}_1\preceq{\mathcal P}_2$ means that the diagram for ${\mathcal P}_2$ is obtained from the diagram for ${\mathcal P}_1$ by adding some (maybe, empty) set of additional dividers. 

Notice also that if $\mathcal P$ is a partition with $m$ terms, then the diagram for $\mathcal P$ contains $m-1$ dividers.\smallskip

If ${\mathcal P}_1\preceq{\mathcal P}_2$, then the {\it interval} $[{\mathcal P}_1,{\mathcal P}_2]$, or the interval with the {\it head} ${\mathcal P}_1$ and the tail ${\mathcal P}_2$, consists of all partition $\mathcal P$ such that ${\mathcal P}_1\preceq{\mathcal P}\preceq{\mathcal P}_2$. In terms of the dots/dividers diagrams, the diagram for a $\mathcal P$ in the interval $[{\mathcal P}_1,{\mathcal P}_2]$ contains all the dividers from the diagram for ${\mathcal P}$ plus some (maybe, empty) subset of the set of the additional dividers of ${\mathcal P}_2$. If $m_i$ is the number of terms of the partition ${\mathcal P}_i$, then the number of these additional dividers is  $(m_2-1)-(m_1-1)=m_2-m_1$, and the set of additional dividers has $2^{m_2-m_1}$ subsets. Accordingly, the interval $[{\mathcal P}_1,{\mathcal P}_2]$ consists of $2^{m_2-m_1}$ partitions. For example, the set of {\it all} ordered partitions of $n$ is the interval with the head $n=n$ and the tail $n=1+1+\dots+1$. So, the total number of ordered partitions of $n$ is $2^{n-1}$.

Notice also that if the number of odd terms in the partition ${\mathcal P}_i$ is $\nu_i$, then the number $\nu$ of odd terms in any partition $\mathcal P$ in the interval $[{\mathcal P}_1,{\mathcal P}_2]$ the inequalities $\nu_1\leq\nu\leq \nu_2$ hold. In particular, if ${\mathcal P}_1$ and ${\mathcal P}_2$ have the same number $\nu$ of odd terms, then every partition in the interval $[{\mathcal P}_1,{\mathcal P}_2]$ also has $\nu$ odd terms.

\subsubsection{Partitions of even and odd types.} 

The following definition looks artificial, but it will be crucially important in the last section of this article. \smallskip

\noindent{\bf Definition.} We say that the ordered partition $n=n_1+n_2+\dots+n_m$ of $n$ belongs to the {\it even (odd) type}, if the first $n_k$, which is not equal to 1, is even (odd). \smallskip

This definition becomes ambiguous for the partition $1+1+\dots+1$, and we declare that this partition belongs to {\it both} even and all types. For example, if $n=5$, then: \smallskip

\noindent Partitions of the even type: $$\begin{array} {c} 4+1,\, 2+3,\ 2+2+1,\  2+1+2,\ 2+1+1+1,\ 1+4,\ 1+2+2,\ 1+2+1+1,\\ 1+1+2+1,\ 1+1+1+2,\ 1+1+1+1+1.\end{array}$$ Partition of the odd type: $$5,\ 3+2,\ 3+1+1,\ 1+3+1,\ 1+1+3,\ 1+1+1+1+1.$$

\subsubsection{Numbers of ordered partitions}

Let $Q_{\rm even}(n), Q_{\rm odd}(n)$ be the numbers of ordered partitions of $n$ of, respectively, even and odd type. Since the full number of ordered partitions of $n$ is $2^{n-1}$ and one partition ($1+1+\dots+1)$ belongs to the both types, the sum $Q_{\rm even}(n)+Q_{\rm odd}(n)$ must be $2^{n-1}+1$. It is easy to find the numbers $Q_{\rm even}(n)$ and $Q_{\rm odd}(n)$ for small values of $n$:\vskip.2in

\centerline{
\beginpicture
\setcoordinatesystem units <1in,1in> point at 0 0
\setplotsymbol({\small.})
\put{$n$} at -1 0
\put{$Q_{\rm even}(n)$} at -1 -.25
\put{$Q_{\rm odd}(n)$} at -1 -.5
\put{$2^{n-1}+1$} at -1 -.75
\plot 1.5 .15 -1.4 .15 -1.4 -.875 1.5 -.875 /
\plot -1.4 -.125 1.5 -.125 /
\plot -1.4 -.375 1.5 -.375 /
\plot -1.4 -.625 1.5 -.625 /
\plot -.6 .15 -.6 -.875 /
\plot -.3 .15 -.3 -.875 /
\plot 0 .15 0 -.875 /
\plot .3 .15 .3 -.875 /
\plot .6 .15 .6 -.875 /
\plot .9 .15 .9 -.875 /
\put1 at -.45 0
\put2 at -.15 0
\put3 at  .15 0
\put4 at .45 0
\put5 at .75 0
\put1 at -.45 -.25
\put2 at -.15 -.25
\put3 at  .15 -.25
\put6 at .45 -.25
\put{11} at .75 -.25
\put1 at -.45 -.5
\put1 at -.15 -.5
\put2 at  .15 -.5
\put3 at .45 -.5
\put6 at .75 -.5
\put2 at -.45 -.75
\put3 at -.15 -.75
\put5 at  .15 -.75
\put9 at .45 -.75
\put{17} at .75 -.75
\put{$\dots$} at 1.1 0
\put{$\dots$} at 1.1 -.25
\put{$\dots$} at 1.1 -.5
\put{$\dots$} at 1.1 -.75
\endpicture
}\vskip.2in

To fill this table in, we can either to list and count all ordered partitions of either  type (at it is done for $n=5$ above), or use the equality $Q_{\rm even}(n)+Q_{\rm odd}(n)=2^{n-1}+1$ and the following simple 

\begin{lemma} For every $n\geq1$, $$Q_{\rm even}(n)=Q_{\rm odd}(n+1).$$\end{lemma}

\noindent{\bf Proof}. To establish a bijection between the even type partitions of $n$ and the odd type partitions of $n+1$, we put $\underbrace{1+1+\dots+1}_n\longleftrightarrow\underbrace{1+1+\dots+1}_{n+1}$ and assign to every other even type partition of $n$ the same partition of $n+1$ with the first even term increased by 1.

Using these equality and lemma, we fill the two middle rows of the table beginning with $Q_{\rm odd}(1)=1\frac\strut\strut$ and following the path \vskip-.16in

\centerline{\hskip2.4in
\beginpicture
\setcoordinatesystem units <1in,1in> point at 0 0
\setplotsymbol(.)
\plot -.45 -.6 -.45 -.35 -.15 -.6 -.15 -.35 .15 -.6 .15 -.35 .45 -.6 .45 -.35 .75 -.6 .75 -.35 /
\arrow <8pt> [.15,.5] from .75 -.35 to 1.05 -.6
\put{\huge.\hskip.06in.\hskip.06in.} at 1.3 -.55
\endpicture
}\vskip.2in

Here is the final result:

\begin{proposition} $$Q_{\rm even}(n)=\left\{\begin{array} {l} \displaystyle{\frac{2^n+2}3},\ {\rm if}\ n\ {\rm is\ even},\\ \\\displaystyle{\frac{2^n+1}3},\ {\rm if}\ n\ {\rm is\ odd}.\;\end{array}\right.\hskip.2in Q_{\rm odd}(n)=\left\{\begin{array} {l} \displaystyle{\frac{2^{n-1}+1}3},\ {\rm if}\ n\ {\rm is\ even},\\ \\ \displaystyle{\frac{2^{n-1}+2}3},\ {\rm if}\ n\ {\rm is\ odd}.\end{array}\right.$$\end{proposition}

\noindent{\bf Proof.} It is sufficient to check that $Q_{\rm odd}(1)=1, Q_{\rm odd}(n)+Q_{\rm even}(n)=2^{n-1}+1$, and $Q_{\rm even}(n)=Q_{\rm odd}(n+1)$. In all these cases the check is immediate.

\begin{corollary} $$Q_{\rm even}(n)=\left\{\begin{array} {ll} 2Q_{\rm odd}(n),& {\rm if}\ n \ {\rm is\ even},\\ 2Q_{\rm odd}(n)+1,& {\rm if}\ n\ {\rm is\ odd}.\end{array}\right.$$\end{corollary} 

    In conclusion, we will provide two more calculations, which will be useful below (in Section 5). It will be the first, but not the last, example of calculations, which result in  the Fibonacci numbers. For the Fibonacci numbers, we use the notation ${\bf F}_n$, where ${\bf F}_0=0, {\bf F}_1=1,$ and ${\bf F}_n-{\bf F}_{n-1}+{\bf F}_{n-2}$, if $n\geq2$. 

We begin with two examples of partitions of even and odd types. First, a partition with all terms equal to 1 or 2 belongs to the even type. Second, a partition with all terms odd belongs to the odd type.

\begin{proposition} The number of ordered partitions of $n$ into 1's and 2's is ${\bf F}_{n+1}$.\end{proposition}

\noindent{\bf Examples.} For $n=1$, $1={\bf F}_2$ partition: 1.

For $n=2$, $2={\bf F}_3$ partitions: $2; 1+1$.

For $n=3$, $3={\bf F}_4$ partitions: $1+2, 2+1;\, 1+1+1$.

For $n=4$, $5={\bf F}_5$ partitions: $2+2; 1+1+2, 1+2+1, 2+1+1;\, 1+1+1+1$.

For $n=5$, $8={\bf F}_6$ partitions: $1+2+2,2+1+2,2+2+1;\, 1+1+1+2, 1+1+2+1,$ $1+2+1+1, 2+1+1+1;\, 1+1+1+1+1$.\smallskip

\noindent{\bf Proof of Proposition 2.6.} Let the number of ordered partitions of $n$ into 1's and 2's be $B_n$. Then, $B_1=1={\bf F}_2,\, B_2=2={\bf F}_3$ (see above). Furthermore, if a partition of $n\geq3$ into 1's and 2's ends with 1, then we remove this 1 and obtain a partition of $n-1$ into 1's and 2's; if it ends with 2, we remove this 2, and obtain a partition  of $n-2$ into 1's and 2's. Moreover, any ordered partition of $n-1$ or of $n-2$ into 1's and 2's may be obtained in this way. Hence, $B_n=B_{n-1}+B_{n-2}$, so $B_n={\bf F}_{n+1}$. 

\begin{proposition} The number of ordered partitions of $n$ into odd terms is $F_{n}$.\end{proposition}

\noindent{\bf Examples.}  For $n=1$, $1={\bf F}_1$ partition: 1.

For $n=2$, $1={\bf F}_2$ partition: $1+1$.

For $n=3$, $2={\bf F}_3$ partitions: $3, 1+1+1$.

For $n=4$, $3={\bf F}_4$ partitions: $1+3,3+1,1+1+1$.

For $n=5$, $5={\bf F}_5$ partitions: $5, 1+1+3, 1+3+1, 3+1+1, 1+1+1+1+1$.

For $n=6$, $8={\bf F}_6$ partitions: $1+5, 3+3, 5+1, 1+1+1+3, 1+1+3+1,$ 

\noindent$1+3+1+1, 3+1+1+1, 1+1+1+1+1+1$.\smallskip

\noindent{\bf Proof of Proposition 2.7.} Let the number of ordered partitions of $n$ into odd terms be $C_n$. Then, $C_1=1={\bf F}_1,\, C_2=1={\bf F}_2$ (see above). Furthermore, if a partition of $n\geq3$ into odd terms ends with 1, then we remove this 1 and obtain a partition of $n-1$ into odd terms; if the last term is 3 or more, we subtract 2 from this term, and obtain a partition  of $n-2$ into odd terms. Moreover, any ordered partition of $n-1$ or of $n-2$ into odd terms  may be obtained in this way. Hence, $C_n=C_{n-1}+C_{n-2}$, so $C_n={\bf F}_n$. 

\subsection{Basic orbits.} 

We will call the orbits corresponding to the sets of non-zero $y$'s {\it basic}. Then, in some sense, all orbits are {\it products of basic orbits.} Let $j_0,\dots,j_m$ and $i_1,\dots,i_m$ be as in Section 2.1. Consider the projection $$G_n\to C_{j_0}\times G_{j_1}\times\dots\times G_{j_m}\eqno(3)$$presented graphically (in the case $n=7,k=2,j_0=2,j_1=3,j_2=2,i_1=2,i_2=5$) in the diagram below:\vskip.1in

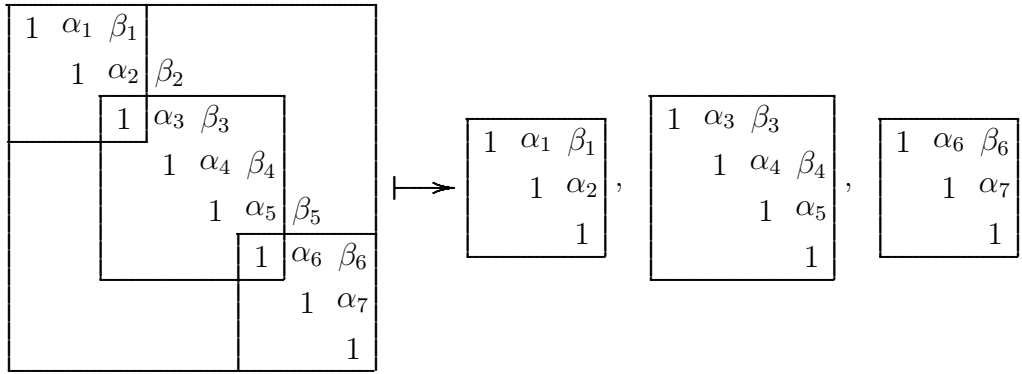
\centerline{
\beginpicture
\setcoordinatesystem units <1.2in,1.2in> point at 0 0
\setplotsymbol({\small.})
\put1 at -2.1 0
\put1 at -1.9  -.2
\put1 at -1.7  -.4
\put1 at -1.5  -.6
\put1 at -1.3  -.8
\put1 at -1.1  -1
\put1 at -.9  -1.2
\put1 at -.7  -1.4
\put {$\alpha_1$} at -1.9 0
\put {$\alpha_2$} at -1.7 -.2
\put {$\alpha_3$} at -1.5 -.4
\put {$\alpha_4$} at -1.3 -.6
\put {$\alpha_5$} at -1.1 -.8
\put {$\alpha_6$} at -.9 -1
\put {$\alpha_7$} at -.7 -1.2
\put {$\beta_1$} at -1.7 0
\put {$\beta_2$} at -1.5 -.2
\put {$\beta_3$} at -1.3 -.4
\put {$\beta_4$} at -1.1 -.6
\put {$\beta_5$} at -.9 -.8
\put {$\beta_6$} at -.7 -1
\plot -1.6  .1 -1.6 -.5 -2.2 -.5 /
\plot -1 -.3 -1.8 -.3 -1.8 -1.1 -1 -1.1 -1 -.3 /
\plot -.6 -.9 -1.2 -.9 -1.2 -1.5 /
\plot -.6 .1 -2.2 .1 -2.2 -1.5 -.6 -1.5 -.6 .1 /
\put 1 at -.1 -.5
\put 1 at .1 -.7
\put 1 at .3 -.9
\put {$\alpha_1$} at .1 -.5
\put {$\alpha_2$} at .3 -.7
\put {$\beta_1$} at ,3 -.5
\plot .4 -.4 -.2 -.4 -.2 -1 .4 -1 .4 -.4 /
\put 1 at .7 -.4
\put 1 at .9 -.6
\put 1 at 1.1 -.8
\put 1 at 1.3 -1
\put {$\alpha_3$} at .9 -.4
\put {$\alpha_4$} at 1.1 -.6
\put {$\alpha_5$} at 1.3 -.8
\put {$\beta_3$} at 1.1 -.4
\put {$\beta_4$} at 1.3 -.6
\plot 1.4 -.3 .6 -.3 .6 -1.1 1.4 -1.1 1.4 -.3 /
\put 1 at 1.7 -.5 
\put 1 at 1.9 -.7 
\put 1 at 2.1 -.9 
\put {$\alpha_6$} at 1.9 -.5
\put {$\alpha_7$} at 2.1 -.7
\put {$\beta_6$} at 2.1 -.5
\plot 2.2 -.4 1.6 -.4 1.6 -1 2.2 -1 2.2 -.4 /
\arrow<8pt> [.15,.5] from -.52 -.7 to -.26 -.7
\plot -.52 -.75 -.52 -.65 /
\put{\title,} at .46 -.7
\put{\title,} at 1.46 -.7
\endpicture
}\vskip.2in

The kernel of the projection (3) is the central subgroup of $G_n$ described by the conditions: all $\alpha$'s and $\beta$'s are zero, except $\beta_{i_1},\dots,\beta_{i_k}$ ($\beta_2$ and $\beta_5$ in the diagram above).

Thus the entries of the group $G_{j_r}$ from this construction ($0\le r\le k$) acquire the notations $\alpha_{i_r+1},\dots,\alpha_{i_{r+1}};\beta_{i_r+1},\dots,\beta_{i_{r+1}-1}$; its adjoint representation is defined in the space of upper triangular matrices with the entries $a_{i_r+1},\dots,a_{i_{r+1}};b_{i_r+1},\dots,b_{i_{r+1}-1}$; and its coadjoint representation is defined in the space of lower triangular matrices with the entries $x_{i_r+1},\dots,x_{i_{r+1}};y_{i_r+1},\dots,y_{i_{r+1}-1}$.  

For fixed non-zero $y_{i_r+1},\dots,y_{i_{r+1}-1}$, plus one additional invariant if $j_r=i_{r+1}-i_r$ is odd (see Section 2.1), there arises a basic orbit of the group $G_{j_r+1}$; fix such orbit for every $r$. The product of these orbits lies the space of the coadjoint representation of the product $G_{j_0+1}\times\dots\times G_{j_k+1}$ and is an orbit of this product group. The projection (3) gives rise to the action of the group $G_n$ in this orbit, and with this action, it is not different from the orbit of $G_n$ corresponding to the partition $n=j_0+\dots+j_k$ and the set of $n-1$ $y$'s, 

\centerline{
\beginpicture
\setcoordinatesystem units <1in,1in> point at 0 0
\setplotsymbol({\small.})
\put{$y_1,\dots,y_{i_1-1},y_{i_1},y_{i_1+1},\dots,y_{i_2-1},y_{i_2},y_{i_2+1},\dots\dots, y_{I_{k-1}-1},y_{i_{k-1}},y_{i_{k-1}+1},\dots,y_{n-1}$} at 0 0
\put{$\|$} at -1.76 -.18
\put{$\|$} at -.4 -.18
\put{$\|$} at 1.4 -.18
\put0 at -1.76 -.38
\put0 at -.4 -.38
\put0 at 1.4 -.38
\endpicture
}
\section{Classes} 

\subsection{Description}

The values of $a_1,\dots,a_n$ are fixed within a class, so we need to consider possible equivalences between the sequences $b_1,\dots,b_{n-1}$ for every fixed sequence $a_1,\dots,a_n$. The sequences $b_1,\dots,b_{n-1}$ and $b'_1,\dots,b'_{n-1}$ are equivalent, if the left system (1),

$$\begin{array} {ccccc} b'_1=b_1&&+\alpha_1a_2&&-\alpha_2a_1\\ b'_2=b_2&&+\alpha_2a_3&&-\alpha_3a_2\\ b'_3=b_3&&+\alpha_3a_4&&-\alpha_4a_3\end{array}$$\vskip-.14in

$$\dots\dots\dots\dots\dots\dots\dots\dots\dots$$\vskip-.18in

$$\begin{array} {l} b'_{n-2}=b_{n-2}+\alpha_{n-2}a_{n-1}-\alpha_{n-1}a_{n-2}\\  b'_{n-1}=b_{n-1}+\alpha_{n-1}a_n\hskip.15in-\alpha_na_{n-1}\end{array}$$

\noindent has a solution in $\alpha_1,\dots\alpha_n$. If all $a$'s are non-zero, then this system has a solution, so in this case there is one class of capacity ${\mathbb K}^{n-1}$ and the set of these classes is labeled by elements of $({\mathbb K^\times})^n$. Let us look, what happens if some $a_i$ are zeroes.

If $a_i=0$, but $a_{i-1}\ne0$ and $a_{i+1}\ne0$ (in particular, $i\ne 1,n$), then the equations with $b_{i-1}$ and $b_i$ become$$\begin{array} {rll} b'_{i-1}=&b_{i-1}&-\alpha_ia_{i-1}\\ b'_i=&b_i&+\alpha_ia_{i+1},\end{array}$$which implies$$b'_{i-1}a_{i+1}+b'_ia_{i-1}=b_{i-1}a_{i+1}+b_ia_{i-1}.$$Thus for the consistency of our system it is necessary that the last equality holds. In other words, $b_{i-1}a_{i+1}+b_ia_{i-1}$ is an invariant, that is, it is fixed within a class.

 We say that the set $\{a_{i+1},a_{i+2},\dots,a_{i+m}\}$ forms a {\it string of zeroes of length $m$} if the following holds: $a_{i+1}=a_{i+2}=\dots=a_{i+m}=0$; if $i>0$ then $a_i\ne0$; if $i+m<n$ then $a_{i+m+1}\ne0$. If $\{a_{i+1},a_{i+2},\dots,a_{i+m}\}$ is string of zeroes of length $m\ge2$, then our system contains a part
 $$\begin{array}{rll} b'_i=&b_i&-\alpha_{i+1}a_i\\  b'_{i+1}=&b_{i+1}\\ b'_{i+2}=&b_{i+2}\end{array}$$ $$\hskip2in\dots\dots\dots\dots\dots\dots\dots\dots\dots\dots$$ $$\hskip.45in\begin{array} {rll} b'_{i+m-1}=&b_{i+m-1}\\ b'_{i+m}=&b_{i+m}&+\alpha_{i+m}a_{i+m+1}\end{array}$$(the first line is absent if $i=0$, and the last line is absent, if $i=n-m$). Thus for the consistency of our system it is necessary that $b'_{i+1}=b_{i+1},\dots,b'_{i+m-1}=b_{i+m-1}$. In other words, $b_{i+1},\dots,b_{i+m-1}$ are invariants.
 
 Below, we will refer to non-zero $a_i$ as to $a$-invariants, and to invariants involving $b_j$ as to $b$-invariants.

 Thus, for any fixed $a_1,a_2,\dots,a_{n-1}$ we have a full description of invariants; this description depends only on the locations of zero $a$'s. For example, if $n=10$ and $a_1=a_3=a_4=a_5=a_7=a_9=a_{10}=0$, while $a_2\ne0,a_6\ne0,a_8\ne0$, then, in the addition to $a$-invariants $a_2,a_6,a_8$,  there are 4 $b$-invariants:  $b_3,b_4,b_6a_8+b_7a_6,b_9$, and a class is determined by fixing their values. Hence in this case classes are parallel 5-dimensional affine planes in the 9-dimensional space spanned by $b_1,b_2,\dots,b_9$.
 
 \subsection{The number of classes of given dimension.}
 
 In Section 2.2, we proved a compact formula for the number of $G_n$-orbits of a given dimension. Unfortunately, no formula of this quality exists for the classes. Below, we restrict ourselves to an inductive procedure of calculating the number of classes of a given dimension. We again assume that $\mathbb{K=F}_q$. An element of ${\mathfrak g}_n$ is characterized by $2n-1$ elements of ${\mathbb F}_q\colon a_1,\dots,a_n,b_1,\dots,b_{n-1}$. The string $a_1,\dots,a_n$ is fixed within any class. We will usually label such string by a sequence of heavy dots $\bullet$ and light dots $\circ$: heavy dots correspond to non-zero $a_i$ and light dots correspond to zero $a_i$. We will denote the number of $k$-dimensional $G_n$-classes by $d_n(k)$. More specifically, $d_n^\bullet(k)$ and $d_n^\circ(k)$ will denote the number of $k$-dimensional  $G_n$-classes with, respectively, $a_1\ne0$ and $a_1=0$. (Thus, $d_n(k)=d_n^\bullet(k)+d_n^\circ(k)$).
 
 For a given $n$, there are $2^n$ strings of heavy and light dots. For each such string we denote the number of heavy dots by $m$ and the number of invariants, calculated as in Section 3.1, by $\ell$. Then this string contributes $(q-1)^mq^\ell$ classes into $d_n^\bullet(n-2-\ell)$ or $d_n^\bullet(n-2-\ell)$, if, respectively, the first dot is heavy or light. For example, the string $$\circ\bullet\circ\circ\bullet\circ\bullet\circ\circ: n=10,m=3,\ell=3\eqno(4)$$contributes $(q-1)^3q^3$  classes into $d_{10}^\circ(5)$.
 
 How to find $d_n(k)$  or, separately, $d_n^\bullet(k)$ and $d_n^\circ(k)$? There is a ``direct" way of doing that: to consider all $2^{n-1}$ strings of heavy and light dots, repeat for each of them the computations similar to those in (4) and calculate the appropriate sums. Let us do this, for example, for $n=3$:
 
 \centerline{
\beginpicture
\setcoordinatesystem units <1in,1in> point at 0 0
\put{$\bullet$} at -1.6 -.2
\put{$\bullet$} at -1.8 -.4
\put{$\bullet$} at -2 -.6
\put{$\bullet$} at -1.8 -.8
\put{$\bullet$} at -1.6 -.8
\put{$\bullet$} at -2 -1
\put{$\bullet$} at -1.6 -1
\put{$\bullet$} at -2 -1.2 
\put{$\bullet$} at -1.8 -1.2
\put{$\bullet$} at -2 -1.4
\put{$\bullet$} at -1.8 -1.4
\put{$\bullet$} at -1.6 -1.4
\put{$\circ$} at -2 0
\put{$\circ$} at -1.8 0
\put{$\circ$} at -1.6 0
\put{$\circ$} at -2 -.2
\put{$\circ$} at -1.8 -.2
\put{$\circ$} at -2 -.4
\put{$\circ$} at -1.6 -.4
\put{$\circ$} at -1.8 -.6
\put{$\circ$} at -1.6 -.6
\put{$\circ$} at -2  -.8
\put{$\circ$} at -1.8 -1
\put{$\circ$} at -1.6 -1.2
\put{$m=0,$} at -1.2 0
\put{$m=1,$} at -1.2 -.2
\put{$m=1,$} at -1.2 -.4
\put{$m=1,$} at -1.2 -.6
\put{$m=2,$} at -1.2 -.8
\put{$m=2,$} at -1.2 -1
\put{$m=2,$} at -1.2 -1.2
\put{$m=3,$} at -1.2 -1.4
\put{$\ell=2,$} at -.7 0
\put{$\ell=1,$} at -.7 -.2
\put{$\ell=0,$} at -.7 -.4
\put{$\ell=1,$} at -.7 -.6
\put{$\ell=0,$} at -.7 -.8
\put{$\ell=1,$} at -.7 -1
\put{$\ell=0,$} at -.7 -1.2
\put{$\ell=0,$} at -.7 -1.4
\put{$q^2$} at -.3 0
\put{$q(q-1)$} at -.1 -.2
\put{$q-1$} at -.2 -.4
\put{$q(q-1)$} at -.1 -.6
\put{$(q-1)^2$} at -.12 -.8
\put{$q(q-1)^2$} at -.08 -1
\put{$(q-1)^2$} at -.1 -1.2
\put{$(q-1)^3$} at -.1 -1.4
\put{in} at .4 0
\put{in} at .4 -.2
\put{in} at .4 -.4
\put{in} at .4 -.6
\put{in} at .4 -.8
\put{in} at .4 -1
\put{in} at .4 -1.2
\put{in} at .4 -1.4
\put{$d_3^\circ(0)$} at .8 0
\put{$d_3^\circ(1)$} at .8 -.2
\put{$d_3^\circ(2)$} at .8 -.4
\put{$d_3^\bullet(1)$} at .8 -.6
\put{$d_3^\circ(2)$} at .8 -.8
\put{$d_3^\bullet(1)$} at .8 -1
\put{$d_3^\bullet(2)$} at .8 -1.2
\put{$d_3^\bullet(2)$} at .8 -1.4
\endpicture
}

From this:$$\hskip-.3in\begin{array} {lll} d_4^\circ(0)=q^2&d_4^\bullet(0)=0&d_4(0)=q^2\\ d_4^\circ(1)=q(q-1)&d_4^\bullet(1)=q(q-1)+q(q-1)^2=q^2(q-1)&d_4(1)=q(q^2-1)\\d_4^\circ(2)=(q-1)+(q-1)^2=q(q-1)&d_4^\bullet(2)=(q-1)^2+(q-1)^3=q(q-1)^2&d_4(2)=q^2(q-1)\end{array}$$However, calculation like this for large $n$ is hardly possible. Below, we deduce an expression of $d_n$ in terms of $d_{n-1}$ and $d_{n-2}$, which can be used for calculating $d_n$ step by step. 

\begin{theorem} For $n\ge3, k\le n-1$,
$$\begin{array} {ll} d_n^\circ(k)&=d_{n-1}^\bullet(k-1)+d_{n-1}^\circ(k)\\ d_n^\bullet(k)&=d_{n-1}^\bullet(k-1)+d_{n-2}^\circ(k-1)+d_{n-2}^\bullet(k-1)\\ &=d_{n-1}^\bullet(k-1)+d_{n-2}(k-1)\end{array}$$
\end{theorem} 

\noindent{\sl Proof\ }  Let us look how $n,m$, and $\ell$ in (4) change, when we cut off one or two first dots from the string of heavy and light dots. \vskip.1in
 
\centerline{
\beginpicture
\setcoordinatesystem units <1in,1in> point at 0 0
\setplotsymbol({\small.})
\plot -2.5 .12 3.8 .12 3.8 -.92 -2.5 -.92 -2.5 .12 /
\plot -1.66 .12 -1.66 -.92 /
\plot -.3 .12 -.3 -.92 /
\plot 1.1 .12 1.1 -.92 /
\plot 2.5 .12 2.5 -.92 /
\put{$n$} at -1 .24
\put{$m$} at .4 .24
\put{$\ell$} at 1.8 .24
\put{$k=n-1-\ell$} at 3.2 .24
\put{$\circ|\bullet\dots$} at -2.18 0
\put{decreases by 1} at -1 0
\put{stays the same} at .4 0
\put{stays the same} at 1.8 0
\put{decreases by 1} at 3.2 0
\put{$\circ|\circ\dots$} at -2.18 -.2
\put{decreases by 1} at -1 -.2
\put{stays the same} at .4 -.2
\put{decreases by 1} at 1.8 -.2
\put{stays the same} at 3.2 -.2
\put{$\bullet|\bullet\dots$} at -2.18 -.4
\put{decreases by 1} at -1 -.4
\put{decreases by 1} at .4 -.4
\put{stays the same} at 1.8 -.4
\put{decreases by 1} at 3.2 -.4
\put{$\bullet\circ|\bullet\dots$} at -2.1 -.6
\put{decreases by 2} at -1 -.6
\put{decreases by 1} at .4 -.6
\put{decreases by 1} at 1.8 -.6
\put{decreases by 1} at 3.2 -.6
\put{$\bullet\circ|\circ\dots$} at -2.1 -.8
\put{decreases by 2} at -1 -.8
\put{decreases by 1} at .4 -.8
\put{decreases by 1} at 1.8 -.8
\put{decreases by 1} at 3.2 -.8
\endpicture
}\vskip.2in

Since every string of heavy and light dots with $n\ge4$ begins with one of the combinations in the first column of the table above, the statement of Theorem follows. 

The relations in Theorem 3.1 determine the induction step. The base of induction is provided by easy to check equalities$$\begin{array} {c} {\rm For}\ n\ge2, d_n^\circ(0)=q^{n-2}, d_n^\bullet(0)=0;\\ d_3^\circ(1)=q-1,d_3^\bullet(1)=q(q-1).\end{array}$$

\subsection{M-classes.} This section, like Section 2.3 above, will not be used before Section 5.

\subsubsection{Sparse sequences and the definition of M-classes.}

A subset $\{i_1,i_2,\dots,i_k\}$ of $\{1,2,\dots,n\}$ is called a sparse sequence, if $i_{s}-i_{s-1}>1$ for every $s, 1<s\leq k$. Let us count the number of sparse sequences.

\begin{proposition} For a given $n$, the number of sparse sequences is the Fibonacci number ${\bf F}_{n+2}$.\end{proposition}

\noindent{\bf Examples.} For $n=1$, there are $2={\bf F}_3$ sparse sequences: $\emptyset,\{1\}$.

For $n=2$, there are $3={\bf F}_4$ sparse sequences: $\emptyset,\{1\},\{2\}$.

For $n=3$, there are $5={\bf F}_5$ sparse sequences: $\emptyset,\{1\},\{2\},\{3\},\{1,3\}$.

For $n=4$, there are $8={\bf F}_6$ sparse sequences: $\emptyset,\{1\},\{2\},\{3\},\{4\},\{1,3\},\{1,4\},\{2,4\}$.\smallskip

\noindent{\bf Proof.} For $n=1,2$ see above. For $n\geq3$, let us consider a sparse sequence in $1,2,\dots,n$. If the last term of this sequence is not $n$, then it is also a sparse sequence in $1,2,\dots,n-1$. If the last term is $n$, then the sequence, being sparse, does not contain $n-1$. So, if we remove the term $n$, we obtain a sparse sequence in $1,2,\dots,n-2$. Thus the number of sparse sequences in $1,2,\dots,n$ is the sum of the numbers of sparse sequences in $1,2,\dots,n-1$ and sparse sequences in $1,2,\dots,n-2$.\smallskip

{\bf Definition.} A class is called an $M$-class, if the subscripts $i_1,i_2,\dots i_k$ of $a$-invariants $a_{i_1},a_{i_2},\dots,a_{i_k}$ form a sparse sequence.

\subsubsection{Containers} For a sparse sequence $I=\{i_1,i_2,\dots,i_k\}\subset\{1,2,\dots,n\}$, we denote by $C(I)$ the subset of $G_n$, which consists of all matrices $g(a_1,\dots,a_n;b_1,\dots,b_{n-1})$ with $a_i\ne0\Longleftrightarrow i\in I$. Thus, $C(I)$ is a union of $M$-classes. Containers play for $M$-classes a role similar to the role which ordered partitions play for orbits: they place together $M$-classes with some important properties being the same. First of all, the $M$-classes from the same container have the same stabilizer.

 The formulas for the multiplication and inversion in the group $G_n$ imply the following formula for the conjugation:$$\begin{array} {l} g(a'_1,\dots,a'_n;b'_1,\dots,b'_{n-1})^{-1}g(a_1,\dots,a_n;b_1,\dots,b_{n-1})g(a'_1,\dots,a'_n;b'_1,\dots,b'_{n-1})\hskip1in\\ \hskip1.6in =g(a_1,\dots,a_n; b_1+a_1a'_2-a_2a'_1,\dots,b_{n-1}+a_{n-1}a'_n-a_na'_{n-1})\end{array}$$Hence $g(a'_1,\dots,a'_n;b'_1,\dots,b'_{n-1})$ belongs to the stabilizer of $g(a_1,\dots,a_n;b_1,\dots,b_{n-1})$ if and only if$$a_1a'_2-a_2a'_1=\dots=a_{n-1}a'_n-a_na'_{n-1}=0,$$which means precisely that, for every $i\in I$, $a'_{i-1}=a'_{i+1}=0$ (here we use the fact that $I$ is sparse and mean that $a'_0=a'_{n+1}=0$). This shows that all elements of $C(I)$ have the same stabilizer, in other words, all $M$-classes within $C(I)$ have the same stabilizer, and this stabilizer $\Stab(I)$ is a normal subgroup of $G_n$. For a more convenient description of $\Stab(I)$, we need new notation. Let $I^-$ be the set of those $i\in\{1,\dots,n\}$, for which $i-1$ or $i+1$ is contained in $I$, and let $I^+$ be the complement of $I^-$ in $\{1,\dots,n\}$. Then $\Stab(I)$ is the set of all $g(a_1,\dots,a_n;b_1,\dots,b_{n-1})$ with $a_i=0$, if $i\in I^-$, that is, $a_i$ may be different from zero only if $i\in I^+$. (Notice that $I^+$ may be not sparse and that $I\subseteq I^+$.)

The group $\Stab(I)$ may be non-commutative. Its commutator subgroup consists of all $g(a_1,\dots,a_n;b_1,\dots,b_{n-1})$ with all $a$'s and $b$'s being zero, with a possible exception of those $b_j$, for which $j,j+1\in I^+$. 

Our description of the container $C(I)$, classes $c\in C(I)$, and the stabilizer $\Stab(I)$ give rise to direct computations of cardinalities of several related sets.

\begin{proposition} {\rm(i)} The number of containers is ${\bf F}_{n+2}$. 

{\rm(ii)} $|\Stab(I)|=q^{n-1+|I^+|}$.

{\rm(iii)} For a class $c\in C(I)$, $|c|=q^{|I^-|}.$

{\rm(iv)} For a class $c\in C(I)$, the number of $b$-invariants is $|I^+|+1$.

{\rm(v)} The number of classes in $C(I)$ is $(q-1)^{|I|}q^{|I^+|-1}$.
\end{proposition}

{\bf Proof.} Part (i) is the same as Proposition 3.2.

Part (ii) follows from our description of $\Stab(I)$: an element of this group is determined by $a_{i_1},\dots,a_{i_k},b_1,\dots,b_{n-1}\in{\mathbb F}_q$. 

Part (iii) follows from Part (ii), since for $c\in C(I),$ $|c|=\displaystyle{\frac{|G_n|}{|\Stab(i)|}=\frac{q^{2n-1}}{|q^{n-1}+|I^|+1}}=q^{n-|I^+|}=q^{|I^-|}.$

$|c|$ can be calculated also directly form the description of $c$: the entries $a_i$ are fixed within $c$, and the entries $b_j$ belong to a subspace of the $b$-space, whose codimension is the number $|b{\rm-inv.}|$  of $b$-invariants. Hence, $|c|=q^{n-1-|b{\rm-inv.}|}$. Compare the two computations of $|c|$: $q^{|I^-|}=q^{n-1-|b{\rm-inv.}|}\Rightarrow 
|I^-|=n-1-|b{\rm-inv.}|\Rightarrow|b{\rm-inv.}|=n-1-|I^-|=|I^+|-1$., as stated in Part (iii).

Part (iv) follows from Part(iii): a class in $C(I)$ is determined by $|I|$ (non-zero) $a$-invariants and $|I^+|-1$ $b$-invariants (which are allowed to be zero).

\section {Representations.} 

The orbit method, as it is described in \cite{kir}, provides a bijective correspondence between coadjoint orbits and unireps.  Our goal is to show that it works for the group $G_n=TD(n,{\mathbb F}_q)$.

 We begin with the observation that the construction of a complex representation corresponding to a given orbit works for this group. We will provide an explicit construction of a representation of $G_n$ corresponding to an orbit and an explicit calculation of the character of this representation. 
 
Then we will show that all the representations constructed are irreducible, and the set of these irreducible representations is complete in the sense that any irreducible complex representation of $G_n$ is equivalent to precisely one of the representations constructed.

Our construction will involve a fixed non-trivial homomorphism ${\bf e}\colon{\mathbb F}_q\to{\mathbb C}^\times$ of the additive group of the field ${\mathbb F}_q$ into the multiplicative group ${\mathbb C}^\times$. Notice that the image of this homomorphism is contained in the set of degree $p$ roots of 1. Notice also that this homomorphism is not unique: there are $q-1$ such homomorphisms labeled by non-zero elements of ${\mathbb F}_q$: the homomorphism ${\bf e}_\alpha\colon{\mathbb F}_q\to{\mathbb C}^\times\ (\alpha\in{\mathbb F}_q-0)$ acts by the formula ${\bf e}_\alpha(\xi)={\bf e}(\alpha\xi).$
 
 \subsection{The construction of representations by the orbit method.} The dimension of a representation corresponding to a $2k$-dimensional orbit is $q^k$. The construction consists of three steps. Let $\Omega$ be an orbit.\smallskip
   
 {\bf First step.} We choose a representative $F\in\Omega$, compute the stabilizer $\Stab F\subset G_n$, and choose a {\it polarization subgroup} $H\subset G_n$ corresponding to some ${\mathfrak h\subset g}_n$ such that $$H\supset\Stab F,\ \dim H=\frac{\dim G_n+\dim\Stab F}2.\ F|_{[{\mathfrak h,h}]}=0.$$The last property shows that the formula $\rho(h)={\bf e}(\tr(Fh))$ defines a 1-dimensional unitary representation of the group $H$.\smallskip
 
{\bf Second step.} Consider the right homogeneous $G_n$-space $T=H\backslash G_n$ and the natural projection $p\colon G_n\to T, p(g)=Hg.$ If we choose for every class $t\in T$ a representative $s(t)\in G_n$, then every element $g\in G_n$ can be uniquely written in the form $g=hs(t)$ where $t\in T,h\in H$. Clearly, here $t$ is the class $Hg\in T$ and $h=gs(t)^{-1}$. Thus the function $s\colon T\to G_n$ allows to identify $G_n$ with the direct product $H\times T$: an element $g\in G_n$ corresponds to the pair $(gs(t)^{-1},t)$. \smallskip

{\bf Third step.} Introduce the so-called {\it Master Equation}: $$s(t)g=h's(t')\ {\rm for\ given}\ t\in T, g\in G_n\ {\rm and\ unknown}\ t'\in T, h'\in H.$$

We have seen above that this equation has a unique solution, $t'=p(s(t)g),h'=s(t)gs(t')^{-1}$.\smallskip

Now we can describe the representation $\pi_\Omega$ of $G_n$. The space of this representation is the space of complex-valued functions on $T$;  since $T$ has dimension $k$ over ${\mathbb F}_q$, the (complex) dimension of the space of function is $q^k$. The representation is described by the formula $$\left(\pi_\Omega(g)f\right)(t)=\rho(h')f(t').$$

\subsection{Basic representations and their characters.} We begin with the representations of $G_n$ corresponding (in the sense of Section 4.1) to the basic orbits (Section 2.4). We begin with the (simpler) case of even $n$.

\subsubsection{The case of even n.} In this case, basic orbits are labeled by sets $(y_1,\dots,y_{n-1})$ of non-zero elements of the field ${\mathbb F}_q$. Values of $x_1,\dots,x_n$ within every orbit are arbitrary, and for a representative $F$ of the orbit corresponding to $(y_1,\dots,y_{n-1})$ we take the matrix\vskip.1in

\centerline{
\beginpicture
\setcoordinatesystem units <1in,1in> point at 0 0
\setplotsymbol({\small.})
\put 0 at -1.4 0
\put 0 at -1.4 -.2
\put{$y_1$} at -1.4 -.4
\put{$y_{n-1}$} at -1 -.8
\put 0 at -.74 -.8
\put 0 at -.54 -.8
\plot -1.32 -.05 -.6 -.73 /
\plot -1.32 -.25 -.8 -.73 /
\plot -1.32 -.45 -1 -.73 /
\plot -1.42 .12 -1.52 .12 -1.52 -.92 -1.42 -.92 /
\plot -.56 -.92 -.46 -.92 -.46 .12 -.56 .12 /
\put{$F=$} at -1.76  -.4
\put{$\Stab F$} at 1 -.22
\put {\title;} at -.36 -.4
\put{consists of} at 1 -.4
\put{matrices} at 1 -.58
\put 1 at 1.6 0
\put 0 at 1.8 0
\put{$\beta_1$} at 2 0
\put 1 at 2.4 -.8
\put 0 at 2.4 -.6
\put{$\beta_{n-1}$} at 2.42 -.4
\plot 1.66 -.08 2.36 -.78 /
\plot 1.84 -.08 2.34 -.58 /
\plot 2.02 -.1 2.26 -.34 /
\plot 1.6 .12 1.5 .12 1.5 -.92 1.6 -.92 /
\plot 2.5 -.92 2.6 -.92 2.6 .12 2.5 .12 /
\endpicture
}\vskip.2in

\noindent (Thus the stabilizer $\Stab F$ of $F$ does not depend on $y_1,\dots, y_{n-1}$.)

For $H$, we can take the group consisting of matrices $h$ shown below (the condition $F|_{[{\mathfrak h,h}]}=0$ holds, because the group $H$ is commutative). The space $T=H\backslash G_n$ has dimension $\displaystyle\frac n2$; the elements of this space are represented by matrices $s(t)=s(t_2,t_4,\dots,t_{n-2},t_n)$ also shown below. \vskip.2in

\centerline{
\beginpicture
\setcoordinatesystem units <1in,1in> point at 0 0
\setplotsymbol({\small.})
\put 1 at -2 0
\put{$\alpha_1$} at -1.8 0
\put{$\beta_1$} at -1.6 0
\put 1 at -1.8 -.2
\put 0 at -1.6 -.2
\put{$\beta_2$} at -1.4 -.2
\put 1 at -1.6 -.4
\put{$\alpha_3$} at -1.4 -.4
\put{$\beta_3$} at -1.2 -.4
\put 1 at -1.4 -.6
\put 0 at -1.2 -.6
\put{$\beta_4$} at -1 -.6
\put{$\ddots$} at -1.2 -.8
\put{$\ddots$} at -.98 -.8
\put{$\ddots$} at -.76 -.8
\put 1 at -1 -1
\put{$\alpha_{n-1}$} at -.76 -1
\put{$\beta_{n-1}$} at -.42 -1
\put 1 at -.7 -1.2
\put 0 at -.4 -1.2
\put 1 at -.4 -1.4
\plot -2 .12 -2.1 .12 -2.1 -1.52 -2 -1.52 /
\plot -.32 -1.52 -.22 -1.52 -.22 .12 -.32 .12 /
\put{$h=$} at -2.3 -.7
\put {\title;} at -.14 -.7
\put{$s(t)=$} at .4 -.7
\plot .82 .12 .7 .12 .7 -1.52 .82 -1.52 /
\put 1 at .8 0
\put 0 at 1 0
\put 0 at 1.2 0
\put 1 at 1 -.2
\put{$t_2$} at 1.2 -.2
\put 0 at 1.4 -.2
\put 1 at 1.2 -.4
\put 0 at 1.4 -.4
\put 0 at 1.6 -.4
\put 1 at 1.4 -.6
\put{$t_4$} at 1.6 -.6
\put 0 at 1.8 -.6
\put{$\ddots$} at 1.6 -.8
\put{$\ddots$} at 1.84 -.8
\put{$\ddots$} at 2.08 -.8
\put 1 at 1.8 -1
\put 0 at 2. -1
\put 0 at 2.26 -1
\put 1 at 2 -1.2
\put{$t_n$} at 2.28 -1.2
\put 1 at 2.28 -1.4
\plot 2.34 -1.52 2.46 -1.52 2.46 .12 2.34 .12 /
\endpicture
}\vskip.2in

The Master equation $s(t)g=h's(t')$ becomes\vskip.2in

\centerline{\hskip.4in
\beginpicture
\setcoordinatesystem units <1in,1in> point at 0 0
\setplotsymbol({\small.})
\put 1 at -2 0
\put 0 at -1.8 0
\put 0 at -1.6 0
\put 0 at -1.4 0
\put 0 at -1.2 0
\put 0 at -2 -.2
\put 1 at -1.8 -.2
\put{$t_2$} at -1.6 -.2
\put 0 at -1.4 -.2
\put 0 at -1.2 -.2
\put 0 at -2 -.4
\put 0 at -1.8 -.4
\put 1 at -1.6 -.4
\put 0 at -1.4 -.4
\put 0 at -1.2 -.4
\put 0 at -2 -.6
\put 0 at -1.8 -.6
\put 0 at -1.6 -.6
\put 1 at -1.4 -.6
\put{$t_4$} at -1.2 -.6
\put 0 at -2 -.8
\put 0 at -1.8 -.8
\put 0 at -1.6 -.8
\put 0 at -1.4 -.8
\put 1 at -1.2 -.8
\put{\bf.} at -2 -1
\put{\bf.} at -1.9 -1
\put{\bf.} at -1.8 -1
\put{\bf.} at -1.7 -1
\put{\bf.} at -1.6 -1
\put{\bf.} at -1.5 -1
\put{\bf.} at -1.4 -1
\put{\bf.} at -1.3 -1
\put{\bf.} at -1.2 -1
\put{\bf.} at -1.1 -1
\put{\bf.} at -1 -1
\put{\bf.} at -1 -.9
\put{\bf.} at -1 -.8
\put{\bf.} at -1 -.7
\put{\bf.} at -1 -.6
\put{\bf.} at -1 -.5
\put{\bf.} at -1 -.4
\put{\bf.} at -1 -.3
\put{\bf.} at -1 -.2
\put{\bf.} at -1 -.1
\put{\bf.} at -1 -.0
\plot -1 .12 -2.12 .12 -2.12 -1 /
\put{\title.} at -.86 -.42
\put 1 at -.6 0
\put{$\alpha_1$} at -.4 0
\put{$\beta_1$} at -.2 0
\put 0 at 0 0
\put 0 at .2 0
\put 0 at -.6 -.2
\put 1 at -.4 -.2
\put{$\alpha_2$} at -.2 -.2
\put{$\beta_2$} at 0 -.2
\put 0 at .2 -.2
\put 0 at -.6 -.4
\put 0 at -.4 -.4
\put 1 at -.2 -.4
\put{$\alpha_3$} at 0 -.4
\put{$\beta_3$} at .2 -.4
\put 0 at -.6 -.6
\put 0 at -.4 -.6
\put 0 at -.2 -.6
\put 1 at 0 -.6
\put{$\alpha_4$} at .2 -.6
\put 0 at -.6 -.8
\put 0 at -.4 -.8
\put 0 at -.2 -.8
\put 0 at 0 -.8
\put 1 at .2 -.8
\plot .34 .12 -.72 .12 -.72 -1 /
\put{\bf.} at -.6 -1
\put{\bf.} at -.5 -1
\put{\bf.} at -.4 -1
\put{\bf.} at -.3 -1
\put{\bf.} at -.2 -1
\put{\bf.} at -.1 -1
\put{\bf.} at 0 -1
\put{\bf.} at .1 -1
\put{\bf.} at .22 -1
\put{\bf.} at .34 -1
\put{\bf.} at .34 -.9
\put{\bf.} at .34 -.8
\put{\bf.} at .34 -.7
\put{\bf.} at .34 -.6
\put{\bf.} at .34 -.5
\put{\bf.} at .34 -.4
\put{\bf.} at .34 -.3
\put{\bf.} at .34 -.2
\put{\bf.} at .34 -.1
\put{\bf.} at .34 -.0
\put{\bf=} at .54 -.42
\plot 1.8 .12 .74 .12 .74 -1 /
\put 1 at .86 0
\put{$\alpha'_1$} at 1.06 0
\put{$\beta'_1$} at 1.26 0
\put 0 at 1.46 0
\put 0 at 1.66 0
\put 0 at .86 -.2
\put 1 at 1.06 -.2
\put 0 at 1.26 -.2
\put{$\beta'_2$} at 1.46 -.2
\put 0 at 1.66 -.2
\put 0 at .86 -.4
\put 0 at 1.06 -.4
\put 1 at 1.26 -.4
\put{$\alpha'_3$} at 1.46 -.4
\put{$\beta'_3$} at 1.66 -.4
\put 0 at .86 -.6
\put 0 at 1.06 -.6
\put 0 at 1.26 -.6
\put 1 at 1.46 -.6
\put 0  at 1.66 -.6
\put 0 at .86 -.8
\put 0 at 1.06 -.8
\put 0 at 1.26 -.8
\put 0 at 1.46 -.8
\put 1 at 1.66 -.8
\put{\bf.} at .86 -1
\put{\bf.} at .96 -1
\put{\bf.} at 1.06 -1
\put{\bf.} at 1.16 -1
\put{\bf.} at 1.26 -1
\put{\bf.} at 1.36 -1
\put{\bf.} at 1.46 -1
\put{\bf.} at  1.56 -1
\put{\bf.} at 1.68 -1
\put{\bf.} at 1.8 -1
\put{\bf.} at 1.8 -.9
\put{\bf.} at 1.8 -.8
\put{\bf.} at 1.8 -.7
\put{\bf.} at 1.8 -.6
\put{\bf.} at 1.8 -.5
\put{\bf.} at 1.8 -.4
\put{\bf.} at 1.8 -.3
\put{\bf.} at 1.8 -.2
\put{\bf.} at 1.8 -.1
\put{\bf.} at 1.8 -.0
\put{\title.} at 1.94 -.42
\put 1 at 2.2 0
\put 0 at 2.4 0
\put 0 at 2.6 0
\put 0 at 2.8 0
\put 0 at 3 0
\put 0 at 2.2 -.2
\put 1 at 2.4 -.2
\put{$t'_2$} at 2.6 -.2
\put 0 at 2.8 -.2
\put 0 at 3 -.2
\put 0 at 2.2 -.4
\put 0 at 2.4 -.4
\put 1 at 2.6 -.4
\put 0 at 2.8 -.4
\put 0 at 3 -.4
\put 0 at 2.2 -.6
\put 0 at 2.4 -.6
\put 0 at 2.6 -.6
\put 1 at 2.8 -.6
\put{$t'_4$} at 3 -.6
\put 0 at 2.2 -.8
\put 0 at 2.4 -.8
\put 0 at 2.6 -.8
\put 0 at 2.8 -.8
\put 1 at 3 -.8
\put{\bf.} at 2.2 -1
\put{\bf.} at 2.3 -1
\put{\bf.} at 2.4 -1
\put{\bf.} at 2.5 -1
\put{\bf.} at 2.6 -1
\put{\bf.} at 2.7 -1
\put{\bf.} at 2.8 -1
\put{\bf.} at 2.9 -1
\put{\bf.} at 3 -1
\put{\bf.} at 3.12 -1
\put{\bf.} at 3.12 -.9
\put{\bf.} at 3.12 -.8
\put{\bf.} at 3.12 -.7
\put{\bf.} at 3.12 -.6
\put{\bf.} at 3.12 -.5
\put{\bf.} at 3.12 -.4
\put{\bf.} at 3.12 -.3
\put{\bf.} at 3.12 -.2
\put{\bf.} at 3.12 -.1
\put{\bf.} at 3.12 0
\plot 3.12 .12 2.08 .12 2.08 -1 /
\put{(5)} at 3.5 -.42
\endpicture
}\vskip.2in

\noindent (the products of matrices are ``truncated": we annihilate all the entries above the $\beta$-diagonal). We equate the entries of the two product matrices and obtain the equalities:$$\begin{array} {c} \alpha'_1=\alpha_1,\, \alpha'_3=\alpha_3,.\dots,\alpha'_{n-1}=\alpha_{n-1};\\ t'_2=t_2+\alpha_2, t'_4=t_4+\alpha_4,\dots,t'_n=t_n+\alpha_n;\\ \beta'_1=\beta_1-t'_2\alpha_1,\beta'_2=\beta_2+t_2\alpha_3,\dots,\beta'_{n-2}=\beta_{n-2}+t_{n-2}\alpha_{n-2},\beta'_{n-1}=\beta_{n-1}-t'_n\alpha_{n-1}\end{array}$$

Now we have a description of the representation $\pi_\Omega$ of $G_n$ corresponding to the basic orbit $\Omega$ corresponding, in turn, to the set $(y_1,\dots,y_{n-1})$ of non-zero elements of ${\mathbb F}_q$. The space of this representation is the $q^{n/2}$-dimensional  space ${\mathcal F}(t_2,t_4,\dots,t_{n-2},t_n)$ of complex-valued functions $f(t_2,t_4,\dots,t_{n-2},t_n)$, and the representation is defined by the formula$$\begin{array} {rl} (\pi_\Omega(g)f)(t_2,t_4,\dots,t_n)&={\bf e}(\tr(Fh'))f(t'_2,t'_4,\dots,t'_n)\\ &={\bf e}(y_1\beta'_1+\dots+y_{n-1}\beta'_{n-1})f(t_2+\alpha_2,t_4+\alpha_4,\dots ,t_n+\alpha_n). \end{array}$$Let us calculate the character of this representation. The space ${\mathcal F}(t_2,t_4,\dots,t_n)$ has a natural basis composed of ``$\delta$-functions": for a fixed set $t_2,t_4,\dots,t_n$, the $\delta$-function takes the value 1 on this set of variables and the value 0 on other sets of variables. The entries of the matrix of the operator $\pi_\Omega(g)$ with respect to this basis correspond to pairs $(t_2,t_4,\dots,t_n),(t'_2,t'_4,\dots,t'_n)$, and this entry is zero, if $(t'_2,\dots,t'_n)\ne(t_2+\alpha_2,\dots,t_n+\alpha_n)$. Thus, if at least one of $\alpha_2,\alpha_4,\dots,\alpha_n$ is not zero, the  matrix of $\pi_\Omega(g)$ has no non-zero diagonal entries, and the trace of this matrix, that is, the value of the character, is zero.

Suppose that $\alpha_2=\alpha_4=\dots=\alpha_n=0$. Then the diagonal entry of the matrix corresponding to $(t_2,\dots,t_n)$ is $$\hskip-.2in\begin{array} {c} {\bf e}(y_1(\beta_1-(t_2+\alpha_2)\alpha_1)+y_2(\beta_2+t_2\alpha_3)+y_3(\beta_3-(t_4+\alpha_4)\alpha_3)+\dots+y_{n-1}(\beta_{n-1}-(t_+\alpha_n)\alpha_{n-1})\\ ={\bf e}(y_1(\beta_1-t_2\alpha_1)+y_2(\beta_2+t_2\alpha_3)+y_3(\beta_3-t_4\alpha_3)+\dots+y_{n-1}(\beta_{n-1}-t_n\alpha_{n-1})=\\ {\bf e}(y_1\beta_1)\dots{\bf e}(y_{n-1}\beta_{n-1}){\bf e}(t_2(y_2\alpha_3-y_1\alpha_1))\dots{\bf e}(t_{n-2}(y_{n-2}\alpha_{n-1}-y_{n-3}\alpha_{n-3})){\bf e}(t_n(-y_{n-1}\alpha_{n-1})).\end{array}$$\vskip.1in

        The sum of all these diagonal entries is $$\begin{array} {l} \displaystyle{{\bf e}(y_1\beta_1)\dots{\bf e}(y_{n-1}\beta_{n-1})\sum_{t_2\in{\mathbb F}_q}{\bf e}(t_2(y_2\alpha_3-y_1\alpha_1))\dots}\\\hskip1in\displaystyle{\sum_{t_{n-2}\in{\mathbb F}_q}{\bf e}(t_{n-2}(y_{n-2}\alpha_{n-1}-y_{n-3}\alpha_{n-3}))\sum_{t_n\in{\mathbb F}_q}{\bf e}(t_n(-y_{n-1}\alpha_{n-1}))}\end{array}$$
        
        Finally,  notice that for a $\gamma\in{\mathbb F}_q$, $$\sum_{t\in{\mathbb F}_q}{\bf e}(t\gamma)=\left\{\begin{array} {l} q,\ {\rm if}\ \gamma=0\\ 0,\ {\rm if}\ \gamma\ne0\end{array}\right.$$ But $y_2\alpha_3-y_1\alpha_1=\dots=y_{n-2}\alpha_{n-1}-y_{n-3}\alpha_{n-3}=-y_{n-1}\alpha_{n-1}=0,$ only if $\alpha_1=\alpha_3=\dots=\alpha_{n-3}=\alpha_{n-1}=0$. We arrive at the final result: the value of the character of the representation $\pi_\Omega$ on the matrix\vskip.1in

\centerline{
\beginpicture
\setcoordinatesystem units <1in,1in> point at 0 0
\setplotsymbol({\small.})
\put 1 at -2 0
\put{$\alpha_1$} at -1.8 0
\put{$\beta_1$} at -1.6 0
\put 1 at -1.8 -.2
\put{$\alpha_2$} at -1.6 -.2
\put{$\beta_2$} at -1.4 -.2
\put{$\ddots$} at -1.6 -.4
\put{$\ddots$} at -1.4 -.4
\put{$\ddots$} at -1.2 -.4
\put 1 at -1.4 -.6
\put{$\alpha_{n-1}$} at -1.16 -.6
\put{$\beta_{n-1}$} at -.82 -.6
\put 1 at -1.1 -.8
\put{$\alpha_n$} at -.8 -.8
\put 1 at -.8 -1
\plot -2 .12 -2.1 .12 -2.1 -1.12 -2 -1.12 /
\plot -.72 -1.12 -.62 -1.12 -.62 .12 -.72 .12 /
\endpicture
}\vskip.2in

\noindent is $q^{n/2}{\bf e}(y_1\beta_1)\dots{\bf e}(y_{n-1}\beta_{n-1}),$ if $\alpha_1=\alpha_2=\dots=\alpha_n=0$, and is zero otherwise.

\subsubsection{The case of odd n.} 

In this case, a basic orbit $\Omega$ is determined by a set $(y_1,\dots,y_{n-1})$ of non-zero elements of the field ${\mathbb F}_q$ {\bf and} a value $I$ of the invariant $$\begin{array} {l}x_1y_2y_4\dots y_{n-1}+y_1x_3y_4y_6\dots y_{n-1}+y_1y_3x_5y_6y_8\dots y_{n-1}+\dots\\ \hskip2in\dots+y_1y_3\dots y_{n-4}x_{n-2}y_{n-1}+y_1y_3\dots y_{n-2}x_n\end{array}\eqno(6)$$For a representative $F$ of $\Omega$ we can take the matrix  \vskip.1in

\centerline{
\beginpicture
\setcoordinatesystem units <1in,1in> point at 0 0
\setplotsymbol({\small.})
\put 0 at -2 0
\put{$x_1$} at -2 -.2
\put{$y_1$} at -2 -.4
\put 0 at -1.8 -.4
\put{$y_{n-1}$} at -1.6 -.8
\put 0 at -1.36 -.8
\put 0 at -1.16 -.8
\plot -1.92 -.5  -1.66 -.7 /
\plot -1.72 -.5  -1.46 -.7 /
\plot -1.92 -.1  -1.22 -.7 /
\plot -2 .1 -2.12 .1 -2.12 -.9 -2 -.9 /
\plot -1.16 .1 -1.04 .1 -1.04 -.9 -1.16 -.9 /
\put{$F=$} at -2.4 -.4
\put{\title,} at -.98 -.4
\put{$x_1=\displaystyle\frac{I}{y_2y_4\dots y_{n-3}y_{n-1}}$} at 0 -.36
\endpicture
}\vskip.2in

The stabilizer $\Stab F$ of $F$ consists of matrices\vskip.2in

\centerline{
\beginpicture
\setcoordinatesystem units <1in,1in> point at 0 0
\setplotsymbol({\small.})
\put{$g(\alpha_1,\dots,\alpha_{n-1};\beta_1,\dots,\beta_{n-2})=$} at -3.3 -.4
\put 1 at -2 0
\put{$\alpha_1$} at -1.8 0
\put{$\beta_1$} at -1.6 0
\put 1 at -1.8 -.2
\put{$\alpha_2$} at -1.6 -.2
\put{$\beta_2$} at -1.4 -.2
\put{$\ddots$} at -1.6 -.4
\put{$\ddots$} at -1.4 -.4
\put{$\ddots$} at -1.2 -.4
\put 1 at -1.4 -.6
\put{$\alpha_{n-1}$} at -1.16 -.6
\put{$\beta_{n-1}$} at -.82 -.6
\put 1 at -1.1 -.8
\put{$\alpha_n$} at -.8 -.8
\put 1 at -.8 -1
\plot -2 .12 -2.1 .12 -2.1 -1.12 -2 -1.12 /
\plot -.72 -1.12 -.62 -1.12 -.62 .12 -.72 .12 /
\endpicture
\hskip.2in}\vskip.2in

\noindent with $\alpha_2=\alpha_4=\dots=\alpha_{n-3}=\alpha_{n-1}=0$ and $\alpha_1,\alpha_3,\dots,\alpha_{n-2},\alpha_n$ satisfying the system of equations $\alpha_1y_1=\alpha_3y_2,\alpha_3y_3=\alpha_5y_4,\dots,\alpha_{n-4}y_{n-4}=\alpha_{n-2}y_{n-3},\alpha_{n-2}y_{n-2}=\alpha_ny_{n-1}$. In other words,$$\begin{array} {rl} \alpha_1&=\kappa\cdot y_2y_4\dots y_{n-3}y_{n-1}\\ \alpha_3&=\kappa\cdot y_1y_4\dots y_{n-3}y_{n-1}\\ \dots\dots&\dots\dots\dots\dots\dots\dot\dots\dots\\ \alpha_{n-2}&=\kappa\cdot y_1y_3\dots y_{n-4}y_{n-1}\\ \alpha_n&=\kappa\cdot y_1y_3\dots y_{n-4}y_{n-2}\end{array}\eqno(7)$$for some $\kappa\in{\mathbb F}_q$. For $H$, we can take the group of matrices shown on the left below (no relations between $\alpha_1,\alpha_3,\dots,\alpha_{n-1}$; again, the condition $F|_{[{\mathfrak h,h}]}=0$ holds, because the group $H$ is commutative). The space $T$ has dimension $\displaystyle\frac{n-1}2$ and the representatives $s(t)=s(t_2,t_4,\dots,t_{n-3},t_{n-1})$ of elements of $T$ are shown below on the right.\vskip.2in

\centerline{
\beginpicture
\setcoordinatesystem units <1in,1in> point at 0 0
\setplotsymbol({\small.})
\put 1 at -2 0
\put{$\alpha_1$} at -1.8 0
\put{$\beta_1$} at -1.6 0
\put 1 at -1.8 -.2
\put 0 at -1.6 -.2
\put{$\beta_2$} at -1.4 -.2
\put 1 at -1.6 -.4
\put{$\alpha_3$} at -1.4 -.4
\put{$\beta_3$} at -1.2 -.4
\put{$\ddots$} at -1.4 -.6
\put{$\ddots$} at -1.2 -.6
\put{$\ddots$} at -1 -.6
\put 1 at -1.2 -.8
\put 0 at -1 -.8
\put{$\beta_{n-1}$} at -.74 -.8
\put 1 at -1 -1
\put{$\alpha_n$} at -.74 -1 
\put 1 at -.78 -1.2
\plot -2 .12 -2.14 .12 -2.14 -1.32 -2 -1.32 /
\plot -.68 .12  -.54 .12 -.54 -1.32 -.68 -1.32 /
\put{\title;} at -.48 -.6
\put 1 at .6 0
\put 0 at .8 0
\put 0 at 1 0
\put 1 at .8 -.2
\put{$t_2$} at 1 -.2
\put 0 at 1.2 -.2
\put 1 at 1 -.4
\put 0 at 1.2 -.4
\put 0 at 1.4 -.4
\put{$\ddots$} at 1.2 -.6
\put{$\ddots$} at 1.4 -.6
\put{$\ddots$} at 1.6 -.6
\put 1 at 1.38 -.8
\put{$t_{n-2}$} at 1.6 -.8
\put 1 at 1.62 -1
\put 0 at 1.82 -1
\put 1 at 1.82 -1.2
\plot .6 .12  .48  .12 .48 -1.32 .6 -1.32 /
\plot 1.82 .12 1.94 .12 1.94 -1.32 1.82 -1.32 /
\put{$s(t)=$} at .2 -.6
\endpicture
}\vskip.2in

The Master equation has the same appearance as before (see (5)), but the matrices now have even, not odd, order. For the entries, we have the following equalities: $$\begin{array} {c} \alpha'_1=\alpha_1,\, \alpha'_3=\alpha_3,.\dots,\alpha'_{n-1}=\alpha_{n-1};\\ t'_2=t_2+\alpha_2, t'_4=t_4+\alpha_4,\dots,t'_{n-1}=t_{n-1}+\alpha_{n-1};\\ \beta'_1=\beta_1-t'_2\alpha_1,\beta'_2=\beta_2+t_2\alpha_3,\dots,\beta'_{n-2}=\beta_{n-2}-t'_{n-1}\alpha'_{n-2},\beta'_{n-1}=\beta_{n-1}+t_{n-1}\alpha_n.\end{array}$$

Now we have a description of the representation $\pi_\Omega$ of $G_n$ corresponding to the basic orbit $\Omega$ corresponding, in turn, to the set $(y_1,\dots,y_{n-1})$ of non-zero elements of ${\mathbb F}_q$ and the value $I\in{\mathbb F}_q$ of the invariant (6). The space of this representation is the $q^{(n-1)/2}$-dimensional  space ${\mathcal F}(t_2,\dots,t_{n-1})$ of complex-valued functions $f(t_2,\dots,t_{n-1})$, and the representation is defined by the formula$$\begin{array} {rl} (\pi_\Omega(g)f)(t_2,\dots,t_{n-1})&={\bf e}(\tr(Fh'))f(t'_2,\dots,t'_{n-1})\\ &={\bf e}(x_1\alpha_1+y_1\beta'_1+\dots+y_{n-1}\beta'_{n-1})f(t_2+\alpha_2,\dots ,t_{n-1}+\alpha_{n-1}). \end{array}$$Let us calculate the character of this representation. The space ${\mathcal F}(t_2,\dots,t_{n-1})$ again has a basis of $\delta$-functions. The entries of the matrix of the operator $\pi_\Omega(g)$ with respect to this basis correspond to pairs $(t_2,\dots,t_{n-2}),(t'_2,\dots,t'_{n-2})$, and this entry is zero, if $(t_2,\dots,t_{n-1})\ne(t_2+\alpha_2,\dots,t_{n-1}+\alpha_{n-1})$. Thus, if at least one of $\alpha_2,\alpha_4,\dots,\alpha_{n-1}$ is not zero, the  matrix of $\pi_\Omega(g)$ has no no-zero diagonal entries, and the trace of this matrix, that is, the value of the character, is zero.

Suppose that $\alpha_2=\alpha_4=\dots=\alpha_n=0$. Then the diagonal entry of the matrix corresponding to $(t_2,\dots,t_n)$ is $$\begin{array} {r} {\bf e}(x_1\alpha_1){\bf e}(y_1(\beta_1-t_2\alpha_1)+y_2(\beta_2+t_2\alpha_3)+\dots\hskip1.6in\\ \dots+y_{n-2}(\beta_{n-2}-t_{n-1}\alpha_{n-2})+y_{n-1}(\beta_{n-1}+t_{n-1}\alpha_n)\end{array}$$ First of all, 
$$\begin{array} {c} {\bf e}(y_1(\beta_1-t_2\alpha_1)+y_2(\beta_2+t_2\alpha_3)+\dots+y_{n-2}(\beta_{n-2}-t_{n-2}\alpha_{n-3})+y_{n-2}(\beta_{n-2}+t_{n-2}\alpha_{n-1})\\ ={\bf e}(y_1\beta_1)\dots{\bf e}(y_{n-1}\beta_{n-1})\cdot  {\bf e}(t_2(-y_1\alpha_1+y_2\alpha_3)\dots{\bf e}(t_{n-1}(-y_{n-2}\alpha_{n-2}+y_{n-1}\alpha_n)\end{array}$$

\noindent The sum of all these diagonal entries is $$\begin{array} {l} \displaystyle{{\bf e}(y_1\beta_1)\dots{\bf e}(y_{n-1}\beta_{n-1})\sum_{t_2\in{\mathbb F}_q}{\bf e}(t_2(y_2\alpha_3-y_1\alpha_1))\dots\sum_{t_{n-1}\in{\mathbb F}_q}{\bf e}(t_{n-1}(-y_{n-2}\alpha_{n-2}+y_{n-1}\alpha_n)}\end{array}$$
 
 Each of the sums in the last formula equals $q$, if $y_2\alpha_3-y_1\alpha_1=\dots=y_{n-1}\alpha_n-y_{n-2}\alpha_{n-2}$, that is, if the equalities (7) hold, and at least one of them is zero otherwise.
 
 It remains to find ${\bf e}(x_1\alpha_1)$ in the case, when (7) holds. But in this case $$x_1\alpha_1=\frac I{y_2\dots y_{n-1}}\cdot\kappa y_2\dots y_{n-1}=\kappa I,\ {\rm so}\ {\bf e}(x_1\alpha_1)={\bf e}(\kappa I).$$
 
 We arrive at the final result: the value of the character of the representation $\pi_\Omega$ on the matrix $g(\alpha_1,\dots,\alpha_n;\beta_1,\dots,\beta_{n-1})$ is $q^{\textstyle\frac{n-1}2}{\bf e}(\kappa I){\bf e}(y_1\beta_1)\dots{\bf e}(y_{n-1}\beta_{n-1}),$ if $\alpha_2=\alpha_4=\dots=\alpha_{n-1}=0$ and the formulas (7) hold for $\alpha_1,\alpha_3,\dots,\alpha_n$. In all other cases, the value of the character is zero.

\subsection{Representations corresponding to all orbits.} We observed in Section 2.4 that all orbits are products of basic orbits. More precisely: if an orbit $\Omega$ of $G_n$ corresponds to an ordered partition $n=j_0+j_1+\dots+j_k$, then this orbit is the products of some basic orbits $\Omega_0,\Omega_1,\dots,\Omega_k$ of the groups $G_{j_0},G_{j_1},\dots,G_{j_k}$ with the action of $G_n$ determined by the projection (3). But in this case the representation $\pi_\Omega$ of $G_n$ determined by these orbits is the tensor product $\pi_{\Omega_0}\otimes\pi_{\Omega_1}\otimes\dots\otimes\pi_{\Omega_k}$ of the basic representations of the groups $G_{j_0},G_{j_1},\dots,G_{j_k}$ corresponding to the orbits $\Omega_0,\Omega_1,\dots,\Omega_k$. And the character of the representation $\pi_\Omega$ is the product of characters of the representations $G_{j_0},G_{j_1},\dots,G_{j_k}$, which were calculated in Section 4.2.\smallskip

Let us prove now that {\bf the representations constructed are all irreducible, and any irreducible representation is equivalent to one of them.}\smallskip

Since all these representations have different characters, they are not equivalent to each other. Let us prove now that

\subsubsection{They are all irreducible.}

Let us begin with basic representations (Section 4.2). The space of each of these representation is the space of complex-valued functions $f(t_2,t_4,\dots,t_n)$ if $n$ is even and $f(t_2,t_4,\dots,t_{n-1})$, if $n$ is odd; the variables $t_i$ in both cases are elements of the field ${\mathbb F}_q$. For a basis of this space we can take the set of ``$\delta$-functions;" each of them takes value 1 for some fixed set of variable $t_i$ and the value 0 on every different set of variables. In is clear that if $f$ is a $\delta$-function concentrated at $(t_2,t_4,t_6,\dots)$, and $g=g(\alpha_1,\dots,\alpha_n;\beta_1,\dots,\beta_{n-1})$, then, $\pi_\Omega(g)f$ is the $\delta$-function concentrated at $(t_2+\alpha_2,t_4+\alpha_4,t_6+\alpha_6,\dots)$ times a non-zero coefficient. This means that if a space of a subrepresentation of our representation contains a $\delta$-function, then this space contains a whole basis, that is, this subrepresentation coincides with the whole representation. 

What if the space of a subrepresentation contains no $\delta$-functions? Take a non-zero element $f$ of this space, which is a linear combination of the minimal possible number $m>1$ of $\delta$-functions (that is, $f$ takes non-zero values at $m$ points). We can apply to it an operator $\pi_\Omega(g)$, where $g$ is a matrix with $\alpha_2=\alpha_4=\alpha_6=\dots=0$, Then the result will be a linear combination of the same $\delta$-function as $f$ (that is, the function $\pi_\Omega(g))f$ takes non-zero values at the same points as $f$). More precisely, the value at $(t_2,t_4,t_6,\dots)$ is multiplied by $${\bf e}(y_1(\beta_1-t_2\alpha_1)+y_2(\beta_2+t_2\alpha_3)+y_3(\beta_3-t_4\alpha_3)+y_4(\beta_4+t_4\alpha_5)+\dots)$$ \noindent ($t'_i=t_i$, since $\alpha_2=\alpha_4=\dots=0$) times ${\bf e}(x_1\alpha_1)$, if $n$ is even. 

Suppose that $f$ takes non-zero values at the points $(t_2,t_4,t_6,\dots)\ne(t'_2,t_4,t'_6,\dots)$.  If $t_2\ne t'_2$, then we apply to this $f$ the matrix $g$ with $\alpha_1\ne0, \alpha_2=\alpha_3=\dots=0$. The values $f(t_2,t_4,t_6,\dots), f(t'_2,t'_4,t'_6,\dots)$ will be multiplied by different numbers. If $t'_2=t_2$, but $t'_4\ne t_4$, then we apply $g$ with $\alpha_3\ne0$ and all other $\alpha$'s are zeroes -- again $f(t_2,t_4,t_6,\dots), f(t'_2,t'_4,t'_6,\dots)$ will be multiplied by different numbers. And so on. Thus $\pi_\Omega(g)f(t_2,t_4,t_6,\dots)=cf(t_2,t_4,t_6,\dots),\pi_\Omega(g)f(t'_2,t'_4,t'_6,\dots)=c'f(t'_2,t'_4,t'_6,\dots),$ and $c'\ne c$. Let $\widetilde f =cf-\pi_\Omega(g)f$. The function $\widetilde f$ is not zero, belongs to the space of the subrepresentation, and  takes non-zero values at less than $m$ points (since $\widetilde f(t_2,t_4,t_6,\dots)=0$). This contradicts to the minimality of $m$. Thus, the space of the subrepresentation contains a $\delta$-function, hence the subrepresentation coincides with the whole representation, hence this whole representation is irreducible. This completes the proof of irreducibility of all basic representations.

The construction of a general representation is described in Section 2.4. We fix an ordered partition$$n=j_1+j_2+\dots+j_k$$of $n$ and consider the projection (3)$$G_n\to G_{j_1}\times G_{j_2}\times\dots\times G_{j_k}$$  (see Section 2.4). Then we fix basic orbits $\Omega_1,\Omega_2,\dots,\Omega_k$ of the groups $G_{j_1}, G_{j_2},\dots, G_{j_k}$ and corresponding representations $\pi_{\Omega_1}, \pi_{\Omega_2},\dots,  \pi_{\Omega_1}$ of these groups. The space $V_i$ of the representation $\pi_{\Omega_i}$ is the space of complex-valued functions of $\displaystyle{\left[\frac{j_i}2\right]}$ variables from ${\mathbb F}_q$. The tensor product $$V=V_1\otimes V_2\otimes\dots\otimes V_k$$ becomes the space of representation $\pi_\Omega$ of $G_{j_1}\times G_{j_2}\times\dots\times G_{j_k}$, and the projection (3)  turns it into a representation of $G_n$; this is a construction of a general representation of $G_n$. Let us prove that this representation is irreducible.

The basis of the space $V_i$ consists of $\delta$-functions. Hence the basis of the space $V$ consists of tensor products$f_1\otimes f_2\otimes\dots\otimes f_k,$where $f_i$ is a delta-function from $V_i$. Take a non-zero subrepresentation of this representation. We want to prove that the space $W$ of this subrepresentation contains a vector from the basis. Take a non-zero vector $f$ of this space. It is a linear combination of the vectors of the basis. Suppose that this linear combination involves more than 1 basic vectors, say, it involves different $f_1\otimes f_2\otimes\dots\otimes f_k$ and $f'_1\otimes f'_2\otimes\dots\otimes f'_k$ (where all $f_i,f'_i$ are $\delta$-functions). Then $f'_i\ne f_i$ for some $i$. We have seen before that there is some $g_i\in G_{j_i}$, which takes every $\delta$-function into itself times a non-zero coefficient and, in particular, $\pi_{\Omega_i}(g_i)f_i=cf_i,\pi_{\Omega_i}(g_i)f'_i=c'f_i,\, c'\ne c$. Let $g=(1,\dots,1,g_i,1\dots,1)\in G_{j_1}\times G_{j_2}\times\dots\times G_{j_k}$. Then $\pi_\Omega(g)f-cf\in W$ is a linear combination of the same basic vectors as $f$ except $f_1\otimes f_2\otimes\dots\otimes f_k$, thus it involves one less basic vectors than $f$. Repeating this operation, we arrive at a basic vector, which is contained in $W$. Appropriate elements of $G_{j_1}\times G_{j_2}\times\dots\times G_{j_k}$ take this basic vector into all other basic vectors. Thus, $W=V$, and our representation is irreducible.

\subsubsection{The list of irreducible representation is complete}

To prove this, we need to check that the sum of the squares of dimension or the representation constructed is equal to the number of elements of the group $G_n$, that is, to $q^{2n-3}$. The dimension of the representation corresponding to a $2k$-dimensional orbit is $q^k$, the number of $2k$-dimensional orbits was calculated in Section 2.2 (Theorem 2.1). Thus,  we need to prove the following equality:$$\sum_{k\ge0}\left(q^{n-k-1}(q-1)^k\left({n-k-1\choose k}q+{n-k-1\choose k-1}\right)\right)\cdot q^{2k}=q^{2n-1}$$or$$\sum_{k\ge0}q^{n+k-1}(q-1)^k\left({n-k-1\choose k}q+{n-k-1\choose k-1}\right)=q^{2n-1}.\eqno(8)$$Put$$P_i=\sum_{k\ge0}q^{i+k}(q-1)^k{i-k-1\choose k}$$(In particular, $P_1=q,P_2=q^2$). Then the left hand side of the equality (8) is the sum of $P_n$ and$$\sum_{k\ge1}q^{n+k-1}(q-1)^k{n-k-1\choose k-1}=\sum_{k\ge0}q^{n+k}(q-1)^{k+1}{n-k-2\choose k}=q(q-1)P_{n-1}.$$On the other hand, $$P_i=\sum_{k\ge0}q^{i+k}(q-1)^k{i-k-1\choose k}=\sum_{k\ge0}q^{i+k}(q-1)^k{i-k-2\choose k}+\sum_{k\ge1}q^{i+k}(q-1)^k{i-k-2\choose k-1}$$ $$=\sum_{k\ge0}q^{i+k}(q-1)^k{i-k-2\choose k}+\sum_{k\ge0}q^{i+k+1}(q-1)^{k+1}{i-k-3\choose k}=qP_{i-1}+q^3(q-1)P_{i-2}.$$

From this: $$P_n+q(q-1)P_{n-1}=(q(q-1)+q)P_{n-1}+q^3(q-1)P_{n-2}=$$ $$q^2P_{n-1}+q^3(q-1)P_{n-2}=(q^3(q-1)+q^3)P_{n-2}+q^5(q-1)P_{n-3}=$$ $$ q^4P_{n-2}+q^5(q-1)P_{n-3}=(q^5(q-1)+q^5)P_{n-3}+q^7(q-1)P_{n-4}=$$  $$q^6P_{n-3}+q^7(q-1)P_{n-4}=\dots=q^{2i}P_{n-i}+q^{2i+1}(q-1)P_{n-i-1}=\dots=$$ $$q^{2n-4}P_2+q^{2n-3}(q-1)P_1=q^{2n-2}+q^{2n-2}(q-1)=q^{2n-1},$$as was stated.

\section{Models}

\subsection{Introduction}

For a finite group $G$, a {\it model} is a finite dimensional representation, whose decomposition into irreducible representations contains an irreducible representation of every equivalence class precisely once. The formal definition of a model (as well as the term {\it model}) was first introduced by I.N.Bernstein, I.M. Gelfand and S. I. Gelfand in their article \cite{bgg}; its significance for the representation theory was demonstrated later by I.M. Gelfand and A.V. Zelevinsky \cite{gz}.

If we have a classification of irreducible representations of $G$ (which is the case for $G=G_n$), then there is a ready construction of a model. What we are interested in here is a ``{\it geometric}" construction of a model for the group $G_n$, the meaning of which we describe below. This description does not use any specific properties of the group $G_n$, so we are speaking of an arbitrary group $G$.

We fix an ad-invariant (that is, consisting of whole classes) set $M\subset G$ and consider the union $\widehat M=\mathop\bigcup\limits_{\gamma\in M}{\mathbb C}_\gamma$ of 1-dimensional complex vector spaces ${\mathbb C}_\gamma\cong{\mathbb C}$. Then we fix a lifting of the adjoint action of $G$ in $M$ to the action in $\widehat M$: for a $g\in G$, the transformation $\widehat g\colon\widehat M\to\widehat M$ maps isomorphically ${\mathbb C}_\gamma$ onto ${\mathbb C}_{g\gamma g^{-1}}$. The space of the representation of $G$, which is the goal of our construction is the space $\mathcal M$ of ``twisted functions" $f\colon M\to\widehat M, f(\gamma)\in{\mathbb C}_\gamma$. The action of $G$ in $\mathcal M$ is defined by formula $gf(\gamma)=\widehat gf(g^{-1}\gamma g).$ Our goal is to choose an $M\subset G_n$ and the lifted action of $G_n$ in $\widehat M$ in such a way that the representation of $G_n$ in $\mathcal M$ be a model for $G_n$.

Return to the general case. There are two ways of understanding construction above: algebraic and geometric.

The algebraic way, which we, actually, follow below, is the following. For every conjugacy class $c\subset M$, choose a representative $\gamma_c\in c$ and consider the stabilizer $\Stab(c)=\{g\in G\vert g\gamma_c=\gamma_cg\}$ of $\gamma_c$. Then for $g\in\Stab(c),\ \widehat g$ maps ${\mathbb C}_{\gamma_c}$ into ${\mathbb C}_{\gamma_c}$, forming a 1-dimensional representation $\rho_c$ of $\Stab(c)$. For $g\in\Stab(c)$, the map ${\mathbb C}_{\gamma_c}\mathop\to\limits^{\widehat g}{\mathbb C}_{\gamma_c}$ is the multiplication by a non-zero complex number $\chi_g$, and the homomorphism $\chi_c\colon\Stab(c)\to{\mathbb C}^\times$ is the character of the representation $\rho_c$. These $\chi_c$ determine (up to an isomorphism) the whole action of $G$ in $\widehat M$, and the representation $\mathcal M$ becomes the direct sum $\mathop\bigoplus\limits_{c\subset M}\Ind_{\Stab(c)}^G\rho_c$ of induced representations.

The geometric way fits better the topological case. The union $\widehat M=\mathop\bigcup\limits_{\gamma\in M}{\mathbb C}_\gamma$ may be regarded as a line bundle over $M$, and the action of $G$ in $\widehat M$ is a fiberwise lifting of the action of $G$ in $M$. Twisted functions are sections of the bundle, and the action is the natural action.

This way of constructing a model was used by A. A. Klyachko \cite{kl} in his construction of a model for the symmetric group $S_n$. For $M\subset S_n$, Klyachko used the set of all ``involutions" $M=\{s\in S_n\vert s^2=1\}$. (Later the significance of the manifold of involutions in the representation theory was demonstrated by Anne Melnikov \cite{Mel}.) In our construction of a model for the group $G_n$, we follow this idea: $M\subset G_n$ is precisely the set of involutions in the case $p=2$. For $p>2$, there are no involutions in $G_n$ (besides the identity), and we need  to modify the definition of $M$. We do it in the section 5.2 below.

\subsection{The case of $\bf G_n$: the set M, stabilizers and characters.}

For the set $M\subset G_n$ (see Section 5.1) we take the set of matrices $g(a_1,\dots,a_n;b_1,\dots,b_{n-1})\in G_n$ with $a_ia_{i+1}=0$ for $i=1,\dots,n-1$. Thus, the subscripts $i$, for which $a_i\neq0$ must form a sparse sequence in $\{1,\dots,n\}$. In other words, $M$ is the union of all $M$-{\it classes}, which we studied in details in Section 3.3. In Section 3.3.2, we formed unions of $M$-classes, which we called {\it containers}. Containers correspond to sparse sequences $I\subset\{1,\dots,n\}$; the container $C(I)$ consists of those $g(a_1,\dots,a_n;b_1,\dots,b_{n-1})$, for which $a_i\ne0$ if and only if $i\in I$. Every $M$-class belongs to one container.The number of containers in $G_n$ is the Fibonacci number $F_{n+2}$.

Recall two useful definitions from Section 3.3. We denote by $I^-$ the set of those $i$, for which either $i+1$ or $i-1$ belongs to $I$. It is obvious that $I^-\cap I=\emptyset.$ The complement $\{1,\dots,n\},-I^-\supset I$ is denoted by $I^+$. A class $c\subset C(I)$ is determined by $|I|$ $a$-invariants $a_i\neq0, i\in I$ and $|I^+|-1$ $b$-invariants, which are allowed to be 0 (see Sections 3.1 and 3.3.2). Thus, the container $C(I)$ contains $(q-1)^{|I|}q^{|I^+|-1}$ classes; every class in $C(I)$ contains $q^{n-|I^+|}$ elements of $G_n$.

 All $g\in C(I)$ have the same stabilizer, which we denote by $\Stab(I)$. We repeat the description of $\Stab(I)$ from Section 3.3.2. Let $I^-$ be the set of those $i$, for which either $i+1$ or $i-1$ belongs to $I$. It is obvious that $I^-\cap I=\emptyset.$ The complement $\{1,\dots,n\},-I^-\supset I$ is denoted by $I^+$. The stabilizer $\Stab(I)$ consists of $g(a_1,\dots,a_n;b_1,\dots,b_{n-1})$ such that $a_i=0$ for $i\in I^-$, that is, $a_i$ may be different from zero only for $i\in I^+$.

The group $\Stab(I)$ is normal, but may be non-commutative. Its commutator subgroup consists of all $g(a_1,\dots,a_n;b_1,\dots,b_{n-1})$ with all $a$'s and $b$'s being zero, with a possible exception of those $b_j$, for which $j,j+1\in I^+$. Hence, every homomorphism $\chi\colon\Stab(I)\to{\mathbb C}^\times$ has the form$$\chi(g(a_1,\dots,a_n;b_1,\dots,b_{n-1}))=\prod{\bf e}(A_ia_i)\cdot\prod{\bf e}(B_jb_j),\eqno(9)$$where $A_i,B_j\in{\mathbb F}_q$, $A_i=0$ for $i\in I^-$, and $B_j=0$, if $j,j+1\in I^+$. (Recall that $\bf e$ is a fixed non-trivial homomorphism of the additive group of the field ${\mathbb F}_q$ into the multiplicative group ${\mathbb C}^\ast$ -- see Section 4.)  Thus,  the construction of a model will involve a choice of the coefficients $A_i$ and $B_j$ for every $M$-class. This choice will be done in Section 5.X.

In conclusion, we provide a formula for the character of the representation of 
$G_n$ induced by a 1-dimensional representation of $\Stab(I)$ withy the character (9).

According to the classical Frobenius formula, If $\chi$ is the character of some finite-dimensional representation of a subgroup $H$ of a finite group $G$, then the character $\widehat\chi$ of the induced representation of $G$ is described by$$\widehat\chi(s)=\frac1{|H|}\sum_{t\in G, t^{-1}st\in H}\chi(t^{-1}st).$$

If the subgroup $H$ is normal (which is our case), then the description becomes $$\widehat\chi(s)=0,\ {\rm if}\ s\notin H,\ {\rm and}\ \widehat\chi(s)=\frac1{|H|}\sum_{t\in G}\chi(t^{-1}st),\ {\rm if}\ s\in H$$

The formulas for the group operations in $G_n$  (Section 1) imply the formula for the conjugation:
$$\begin{array} {l} g(\alpha_1,\dots,\alpha_n; \beta_1,\dots,\beta_{n-1})^{-1}g(a_1,\dots,a_n; b_1,\dots,b_{n-1})g(\alpha_1,\dots,\alpha_n; \beta_1,\dots,\beta_{n-1})\\ \hskip1.6in=g(a_1,\dots,a_n; b_1+a_1\alpha_2-a_2\alpha_1,\dots,b_{n-1}+a_{n-1}\alpha_n-a_n\alpha_{n-1}).\end{array}$$
In our case, the Frobenius formula gives (we use the notations $\overline a=\{a_1,\dots,a_n\}$ and $\overline b=\{b_1,\dots,b_{n-1}\}$:$$\widehat\chi(g(\overline a,\overline b))=\frac1{|\Stab(I)|}\sum_{\alpha_i,\beta_j}\left(\prod{\bf e}(A_ia_i)\cdot\prod{\bf e}(B_j(b_j+a_j\alpha_{j+1}-a_{j+1}\alpha_j))\right).$$ The expression under $\sum_{\alpha_i,\beta_j}$ does not depend on $\beta$'s and on $\alpha_i$ with $i\in I^+$, so the formula may be simplified:$$\widehat\chi(g(\overline a,\overline b))=\chi(g(\overline a,\overline b))\cdot\prod_{j\in I^-}\sum_{\alpha_jj}{\bf e}((B_{j-1}a_{j-1}-B_ja_{j+1})\alpha_j).$$Since $\displaystyle{\sum_y{\bf e}(xy)=\left\{\begin{array} {ll} q,&{\rm if}\ x=0,\\ 0,&{\rm if}\ x\neq0\end{array}\right.}$,we arrive at the final result: $$\widehat\chi(g(\overline a,\overline b))=\left\{\begin{array} {ll} q^{|I^-|}\chi(g(\overline a,\overline b)),&{\rm if\ for\ every}\ j\in I^-,a_j=0\ {\rm and}\  B_{j-1}a_{j-1}=B_ja_{j+1},\\ 0&{\rm otherwise}.\end{array}\right.$$

\subsection{Flocks.}

\subsubsection{Introduction.}

The construction of a model for the group $G_n$, which we outlined in Section 5.1, is supposed to make every irreducible representation of $G_n$ an irreducible component of a certain representation of $G_n$ induced from a 1-dimensional representation of the stabilizer $\Stab(c)$ of a certain M-class $c\subset G_n$. This would determine a mapping \vskip.2in

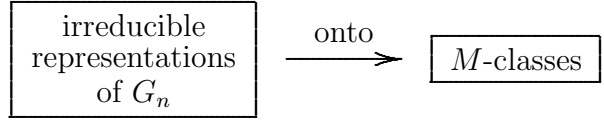
\centerline{
\beginpicture
\setcoordinatesystem units <1in,1in> point at 0 0
\setplotsymbol({\small.})
\put{irreducible} at -1 .18
\put{representations} at -1 0
\put{of $G_n$} at -1 -.18
\put{$M$-classes} at 1 0
\put{onto} at .1 .14
\plot -1.64 .3 -.36 .3 -.36 -.3 -1.64 -.3 -1.64 .3 /
\plot .56 .12 1.44 .12 1.44  -.12 .56 -.12 .56 .12 /
\arrow <10pt> [.15,.5] from -.2 0 to .4 0
\endpicture
}\vskip.2in

This map, actually, will take irreducible representations related to one ordered partition of $n$ (see Section 4.3) into $M$-classes from one container (we will see in Section 5.X that there will be one small exception to this rule), the mapping (10) will be a refinement of a certain mapping\vskip.2in

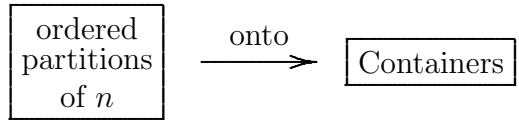
\centerline{
\beginpicture
\setcoordinatesystem units <1in,1in> point at 0 0
\setplotsymbol({\small.})
\put{ordered} at -1 .18
\put{partitions} at -1 0
\put{of $n$} at -1 -.18
\put{Containers} at .8 0
\put{onto} at -.1 .14
\plot -1.4 .3 -.6 .3 -.6 -.3 -1.4 -.3 -1.4 .3 /
\plot .36 .12 1.24 .12 1.24 -.12 .36 -.12 .36 .12 /
\arrow <10pt> [.15,.5] from -.4 0 to .2 0
\endpicture
}\vskip.2in

For every $n$, the number of containers is the Fibonacci number $F_{n+2}$, while the number of partitions is $2^{n-1}$. If $n$ is large enough ($n>4$), then $2^{n-1}>F_{n+2}$ (for example, for $n=6$ these two numbers are 32 and 21, and for $n=10$ they are 512 and 144). So, we must be prepared to the fact that several partitions will be assigned to the same container. For the sets of partitions, which are going to be assigned to one container, we will use the term {\it flocks}. We will begin with a descriptions of flocks. \smallskip

\subsubsection{Description of flocks.} In Section 2.3.2, we introduced a definition of two {\it types} of (ordered)  partitiions: {\it even} and {odd}. Remind that an ordered partition $n=n_1+n_2+\dots+n_m$ belongs to the {\it even (odd) type}, if the first $n_i\neq1$ (if such $n_i$ exists) is even (odd). According to this definition, the partition $n=1+1+\dots+1$ belongs to the both types. 

Remind some useful notations from Section 2.3. The numbers of even and odd terms in the partition of $n$ are denoted, respectively, by $\mu$ and $\nu$; the numbers $n$ and $\nu$ are of the same parity, and number $(n-\nu)/2$ is denoted by $k$.

All the partitions within a flock will be of the same, even or odd, type. Accordingly, we will speak of {\it flocks} of {\it even} or {\it odd} type.  Every flock will be an {\it interval} with respect to the partial ordering of partitions (see Section 2.3.1), that is, it has a {\it head} $\mathcal H$ and a {\it tail} $\mathcal T$ and consists of all partitions $\mathcal P$ such that $\mathcal{H\preceq P\preceq T}$. 

The head of a partition of odd type must be a partition into $\nu=n-2k$ odd parts and no even parts. The tail of a partition of odd type must have the first term greater than 1 (if there is any) equal to 3, and all the other terms equal to 1 or 2. To obtain the tail from the head $n=n_1+n_2+\dots+n_\nu$ we replace the first $n_i\leq3$ by $3+2+\dots+2$ and every other $n_i$ by $1+2+\dots+2$. For example, if the head is $1+7+3+1+3$, then the tail will be $1+\underbrace{3+2+2}+\underbrace{1+2}+1+\underbrace{1+2}$. Every partition $n=n_1+\dots+n_m$ of the odd type belongs to precisely one flock of the odd type. To obtain the head of this flock we need to combine every odd $n_i$ with all even terms after it (and before the next odd term); to obtain the tail of this partition, we need to replace the first odd term greater than 1 by $3+2+\dots+2$, every other odd term by $1+2+\dots+2$ and every even term by $2+2+\dots+2$.

The head of a partition of even type must begin with $1+1+\dots+1+$ even term and have all the terms after that odd. The tail of a partition of the even type must consist of 1's and 2's. To obtain the tail from the head, we replace the even term (if there is any) by $2+2+\dots+2$ and replace every odd term by $1+2+\dots+2$. For example, if the head is $1+1+4+3+5$, then the tail will be $1+1+\underbrace{2+2}+\underbrace{1+2}+\underbrace{1+2+2}$. Every partition $n=n_1+\dots+n_m$ of the even type belongs to precisely one flock of the even type. To obtain the tail of this flock, we need to replace every odd term by $1+2+\dots+2$ and every even term by $2+2+\dots+2$. To obtain the head of this flock we keep the 1's before the first term greater than 1 unchanged, then combine the even terms after these 1's and before the next odd term, and then combine every odd term with all the even terms after it (and before the next odd term. For example, the partition $1+2+4+5+2+3$ of the even type belongs to the flock of the even type with the head $1+6+7+3$ and the tail $1+2+\underbrace{2+2}+\underbrace{1+2+2}+2+\underbrace{1+2}$.

Notice that the number $\nu$ of odd terms (but not the number $\mu$ of even terms!) and hence the number $k$ are fixed within the flock (of any type).

Notice also that it is convenient to count flocks of the odd type by heads and flocks of the even type by tails. Namely, the number of flocks of the odd type is the number of partitions of $n$ into odd parts, which is $F_n$ (Proposition 2.7). And the number of flocks on the even type is the number of partitions of $n$ into 1's and 2's, which is $F_{n+1}$ (Proposition 2.6). Thus the whole number of flocks is $F_n+F_{n+1}=F_{n+2}$, which conveniently coincides with the number of $M$-classes (see Section 3.3.1).

\subsubsection{Pictures and examples} It is convenient to present this procedure on a picture. 

First we draw $n$ heavy dots in line. Then we put dividers between some dots, so the line is divided into parts, and this partition is the head of a flock. To obtain the tail, we should add some dividers, and we draw them as dotted lines. The description of flocks given above provides an instruction for choosing places for dotted dividers. Namely, for flocks  of the even type, we divide every even part into $2+2+\dots+2$ and every odd part into $1+2+\dots+2$. For the flocks of the odd type, the rules are the same with one exception: we divide the leftmost (odd) part of length $\ge3$ into $3+2+\dots+2$. It easy to prove (we leave this to the reader) that the number of dotted lines will be $k-1$. All the partitions in our flock are obtained by using all solid dividers and  some subset of the set of dotted dividers. Since the set of $k-1$ dotted dividers have $2^{k-1}$ subset, our flock will contain $2^{k-1}$ partitions.  

Two examples of this construction are shown below ($n=12$, odd type, and $k=4$ in the picture on the left, and $n=11$, even type, and $k=4$ in the picture on the right.\vskip.17in

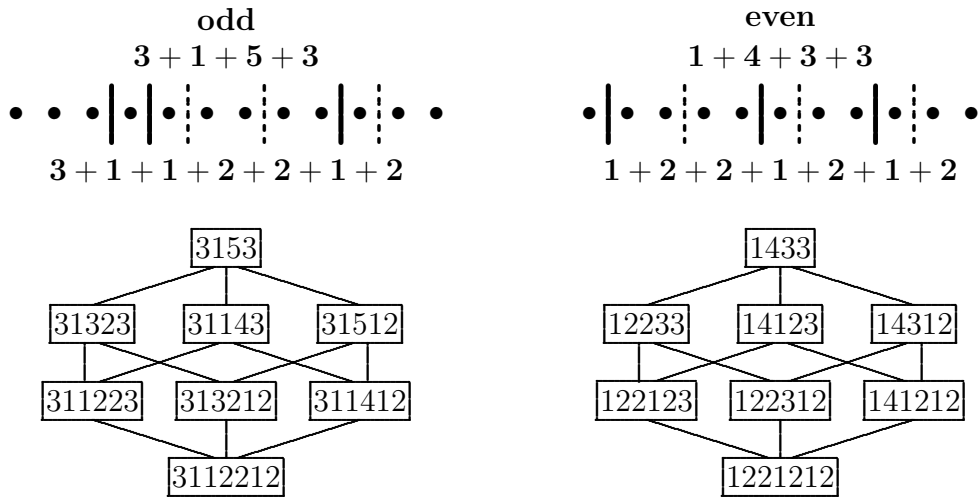
\centerline{
\beginpicture
\setcoordinatesystem units <1in,1in> point   at 0 0
\put{\bf odd} at -1 .2
\put{$\bf3+1+5+3$} at -1 0
\put{$\bullet$} at -2.1 -.3
\put{$\bullet$} at -1.9 -.3
\put{$\bullet$} at -1.7 -.3
\put{$\bullet$} at -1.5 -.3
\put{$\bullet$} at -1.3 -.3     
\put{$\bullet$} at -1.1 -.3
\put{$\bullet$} at -.9 -.3
\put{$\bullet$} at -.7 -.3
\put{$\bullet$} at -.5 -.3
\put{$\bullet$} at -.3 -.3
\put{$\bullet$} at -.1 -.3
\put{$\bullet$} at .1 -.3
\setplotsymbol({\bf.})
\plot -1.6 -.15 -1.6 -.45 /
\plot -1.4 -.15 -1.4 -.45 /
\plot -.4 -.15 -.4 -.45 /
\setplotsymbol(.)
\setdashpattern<2pt,3pt>
\plot -1.2 -.15 -1.2 -.45 /
\plot -.8 -.15 -.8 -.45 /
\plot -.2 -.15 -.2 -.45 /
\put{$\bf3+1+1+2+2+1+2$} at -1 -.6
\setsolid
\setplotsymbol({\small.})
\put{3153} at -1 -1
\plot -1.18 -.9 -.82 -.9 -.82 -1.1 -1.18 -1.1 -1.18 -.9 /
\put{31323} at -1.7 -1.4
\plot -1.92 -1.3 -1.48 -1.3 -1.48 -1.5 -1.92 -1.5 -1.92 -1.3 /
\put{31143} at -1 -1.4
\plot -1.22 -1.3 -.78 -1.3 -.78 -1.5 -1.22 -1.5 -1.22 -1.3 /
\put{31512} at -.3 -1.4
\plot -.52 -1.3 -.08 -1.3 -.08 -1.5 -.52 -1.5 -.52 -1.3 /
\put{311223} at -1.7 -1.8
\plot -1.96 -1.7 -1.44 -1.7 -1.44 -1.9 -1.96 -1.9 -1.96 -1.7 /
\put{313212} at -1 -1.8
\plot -.56 -1.7 -.04 -1.7 -.04 -1.9 -.56 -1.9 -.56 -1.7 /
\put{311412} at -.3 -1.8
\plot -1.26 -1.7 -.74 -1.7 -.74 -1.9 -1.26 -1.9 -1.26 -1.7 /
\put{3112212} at -1 -2.2
\plot -1.3 -2.1 -.7 -2.1 -.7 -2.3 -1.3 -2.3 -1.3 -2.1 /
\plot -1.06 -1.1 -1.7 -1.3 /
\plot -1 -1.1 -1 -1.3 /
\plot -.94 -1.1 -.3 -1.3 /
\plot -1.74 -1.5 -1.74 -1.7 /
\plot -1.66 -1.5 -1.04 -1.7 /
\plot -1.04 -1.5 -1.66 -1.7 /
\plot -.34 -1.5 -.96 -1.7 /
\plot -.96 -1.5 -.34 -1.7 /
\plot -.26 -1.5 -.26 -1.7 /
\plot -1.7 -1.9 -1.06 -2.1 /
\plot -1 -1.9 -1 -2.1 /
\plot -.3 -1.9 -.94 -2.1 /
\put{\bf even} at 1.9 .2
\put{$\bf1+4+3+3$} at 1.9 0
\put{$\bullet$} at .9 -.3
\put{$\bullet$} at 1.1 -.3
\put{$\bullet$} at 1.3 -.3
\put{$\bullet$} at 1.5 -.3
\put{$\bullet$} at 1.7 -.3
\put{$\bullet$} at 1.9 -.3
\put{$\bullet$} at 2.1 -.3
\put{$\bullet$} at 2.3 -.3
\put{$\bullet$} at 2.5 -.3
\put{$\bullet$} at 2.7 -.3
\put{$\bullet$} at 2.9 -.3
\setplotsymbol({\bf.})
\plot 1 -.15 1 -.45 /
\plot 1.8 -.15 1.8 -.45 /
\plot 2.4 -.15 2.4 -.45 /
\setplotsymbol(.)
\setdashpattern<2pt,3pt>
\plot  1.4 -.15 1.4 -.45 /
\plot 2 -.15 2 -.45 /
\plot 2.6 -.15 2.6 -.45 /
\setsolid
\put{$\bf1+2+2+1+2+1+2$} at 1.9 -.6
\setplotsymbol({\small.})
\put{1433} at 1.9 -1
\plot 1.72 -.9 2.08 -.9 2.08 -1.1 1.72 -1.1 1.72 -.9 /
\put{12233} at 1.2 -1.4
\plot .98 -1.3 1.42 -1.3 1.42 -1.5 .98 -1.5 .98 -1.3 /
\put{14123} at 1.9 -1.4
\plot 1.68 -1.3 2.12 -1.3 2.12 -1.5 1.68 -1.5 1.68 -1.3 /
\put{14312} at 2.6 -1.4
\plot 2.38 -1.3 2.82 -1.3 2.82 -1.5 2.38 -1.5 2.38 -1.3 /
\put{122123} at 1.2 -1.8
\plot .94 -1.7 1.46 -1.7 1.46 -1.9 .94 -1.9 .94 -1.7 /
\put{122312} at 1.9 -1.8
\plot 1.64 -1.7 2.16 -1.7 2.16 -1.9 1.64 -1.9 1.64 -1.7 /
\put{141212} at 2.6 -1.8
\plot 2.34 -1.7 2.86 -1.7 2.86 -1.9 2.34 -1.9 2.34 -1.7 /
\put{1221212} at 1.9 -2.2
\plot 1.6 -2.1 2.2 -2.1 2.2 -2.3 1.6 -2.3 1.6 -2.1 /
\plot 1.84 -1.1 1.2 -1.3 /
\plot 1.9 -1.1 1.9 -1.3 /
\plot 1.96 -1.1 2.6 -1.3 /
\plot 1.16 -1.5 1.16 -1.7 /
\plot 1.24 -1.5 1.86 -1.7 /
\plot 1.86 -1.5 1.24 -1.7 /
\plot 2.56 -1.5 1.94 -1.7 /
\plot 1.94 -1.5 2.56 -1.7 /
\plot 2.64 -1.5 2.64 -1.7 /
\plot 1.2 -1.9 1.84 -2.1 /
\plot 1.9 -1.9 1.9 -2.1 /
\plot 2.6 -1.9 1.96 -2.1 /
\endpicture
}\vskip.2in

Below, we show this splitting of the set of partition into flocks for $n=6$ and $n=7$.

\subsubsection{The cases n = 6 and 7} 

\centerline{\bf n  = 6, odd type}
\centerline{
\beginpicture
\setcoordinatesystem units <1in,1in> point   at 0 0
\setplotsymbol({\small.})
\put{$\bf k=2$} at -1.4 .3
\put{15} at -2.2 0
\plot -2.31 .1 -2.09 .1 -2.09 -.1 -2.31 -.1 -2.31 .1 / 
\put{132} at -2.2 -.4
\plot -2.34 -.3 -2.06 -.3 -2.06 -.5 -2.34 -.5 -2.34 -.3 / 
\plot -2.2 -.1 -2.2 -.3 /
\put{$\bf C(1,4,6)$} at -2.2 -.7
\put{33} at -1.4 0
\plot -1.51 .1 -1.29 .1 -1.29 -.1 -1.51 -.1 -1.51 .1 / 
\put{312} at -1.4 -.4
\plot -1.54 -.3 -1.26 -.3 -1.26 -.5  -1.54 -.5 -1.54 -.3 / 
\plot -1.4 -.1 -1.4 -.3 /
\put{$\bf C(1,3,6)$} at -1.4 -.7
\put{51} at -.6 0
\plot -.71 .1 -.49 .1 -.49 -.1 -.71 -.1 -.71 .1 / 
\put{321} at -.6 -.4
\plot -.74 -.3 -.46 -.3 -.46 -.5  -.74 -.5 -.74 -.3 / 
\plot -.6 -.1 -.6 -.3 /
\put{$\bf C(1,3,5)$} at -.6 -.7
\plot -.16 .3 -.16 -.8 /
\put{$\bf k=1$} at 1.1 .2
\put{3111} at .2 -.1
\plot .03 0 .37 0 .37 -.2 .03 -.2 .03 0 /
\put{$\bf C(1,3)$} at .2 -.4
\put{1311} at .8 -.1
\plot .63 0 .97 0 .97 -.2 .63 -.2 .63 0 /
\put{$\bf C(1,4)$} at .8 -.4
\put{1131} at 1.4 -.1
\plot 1.23 0 1.57 0 1.57 -.2 1.23 -.2 1.23 0 /
\put{$\bf C(1,5)$} at 1.4 -.4
\put{1113} at 2 -.1
\plot 1.83 0 2.17 0 2.17 -.2 1.83 -.2 1.83 0 /
\put{$\bf C(1,6)$} at 2 -.4
\plot 2.4 .3 2.4 -.8 /
\put{$\bf k=0$} at 2.8 .2
\put{111111} at 2.8 -.1
\plot 2.54 0 3.06 0 3.06 -.2 2.54 -.2 2.54 0 /
\put{$\bf C(1)$} at 2.8 -.4
\endpicture
}\vskip.2in

\centerline{\bf n  = 6, even type}
\centerline{
\beginpicture
\setcoordinatesystem units <1in,1in> point   at 0 0
\setplotsymbol({\small.})
\put{$\bf k=3$} at -2 .3
\put6 at -2 0
\plot -2.08 .1 -1.92 .1 -1.92 -.1 -2.08 -.1 -2.08 .1 /
\put{24} at -2.16 -.35
\plot -2.27 -.25 -2.05 -.25 -2.05 -.45 -2.27 -.45 -2.27 -.25 /
\put{42} at -1.84 -.35
\plot -1.95 -.25 -1.73 -.25 -1.73 -.45 -1.95 -.45 -1.95 -.25 /
\put{222} at -2 -.7
\plot -2.14 -.6 -1.86 -.6 -1.86 -.8 -2.14 -.8 -2.14 -.6 /
\plot -2.03 -.1 -2.16 -.25 /
\plot -1.97 -.1 -1.84 -.25 /
\plot -2.04 -.6 -2.16 -.45 /
\plot -1.96 -.6 -1.84 -.45 /
\plot -1.5 .4 -1.5 -1.1 /
\put{$\bf C(2,4,6)$} at -2 -1
\put{114} at -1.1 -.1
\plot -1.24 0 -.96 0 -.96 -.2  -1.24 -.2 -1.24 0 /
\put{1122} at -1.1 -.5
\plot -1.27 -.4 -.91 -.4 -.91 -.6 -1.27 -.6 -1.27 -.4 /
\plot -1.1 -.2 -1.1 -.4 /
\put{$\bf C(4,6)$} at -1.1 -.9
\put{$\bf k=2$} at .65 .2
\put{141} at -.4 -.1
\plot -.54 0 -.26 0 -.26 -.2  -.54 -.2 -.54 0 /
\put{1221} at -.4 -.5
\plot -.57 -.4 -.21 -.4 -.21 -.6 -.57 -.6 -.57 -.4 /
\plot -.4 -.2 -.4 -.4 /
\put{$\bf C(3,5)$} at -.4 -.9
\put{411} at .3 -.1
\plot .16 0 .44 0 .44 -.2  .16 -.2 .16 0 /
\put{2211} at .3 -.5
\plot .11 -.4 .49 -.4 .49 -.6 .11 -.6 .11 -.4 /
\plot .3 -.2 .3 -.4 /
\put{$\bf C(2,4)$} at .3 -.9
\put{123} at 1 -.1
\plot .86 0 1.14 0 1.14 -.2  .86 -.2 .86 0 /
\put{1212} at 1 -.5
\plot .83 -.4 1.19 -.4 1.19 -.6 .83 -.6 .83 -.4 /
\plot 1 -.2 1 -.4 /
\put{$\bf C(3,6)$} at 1 -.9
\put{213} at 1.7 -.1
\plot 1.56 0 1.84 0 1.84 -.2  1.56 -.2 1.56 0 /
\put{2112} at 1.7 -.5
\plot 1.51 -.4 1.89 -.4 1.89 -.6 1.51 -.6 1.51 -.4 /
\plot 1.7 -.2 1.7 -.4 /
\put{$\bf C(2,6)$} at 1.7 -.9
\put{231} at 2.4 -.1
\plot 2.26 0 2.54 0 2.54 -.2  2.26 -.2 2.26 0 /
\put{2121} at 2.4 -.5
\plot 2.21 -.4 2.59 -.4 2.59 -.6 2.21 -.6 2.21 -.4 /
\plot 2.4 -.2 2.4 -.4 /
\put{$\bf C(2,5)$} at 2.4 -.9
\endpicture
}
\centerline{
\beginpicture
\setcoordinatesystem units <1in,1in> point   at 0 0
\setplotsymbol({\small.})
\put{$\bf k=1$} at -.4 -1.3
\put{21111} at -2 -1.6
\plot -2.22 -1.5 -1.78 -1.5 -1.78 -1.7 -2.22 -1.7 -2.22 -1.5 /
\put{$\bf C(2)$} at -2 -1.9
\put{12111} at -1.2 -1.6
\plot -1.42 -1.5 -.98 -1.5 -.98 -1.7 -1.42 -1.7 -1.42 -1.5 /
\put{$\bf C(3)$} at -1.2 -1.9
\put{11211} at -.4 -1.6
\plot -.62 -1.5 -.18 -1.5 -.18 -1.7 -.62 -1.7 -.62 -1.5 /
\put{$\bf C(4)$} at -.4 -1.9
\put{11121} at .4 -1.6
\plot .18 -1.5 .62 -1.5 .62 -1.7 .18 -1.7 .18 -1.5 /
\put{$\bf C(5)$} at .4 -1.9
\put{11112} at 1.2 -1.6
\plot .98 -1.5 1.42 -1.5 1.42 -1.7 .98 -1.7 .98 -1.5 /
\put{$\bf C(6)$} at 1.2 -1.9
\plot 1.7 -1.2 1.7 -2 /
\put{$\bf k=0$} at 2.2 -1.3
\put{111111} at 2.2 -1.6
\plot 1.94 -1.5 2.46 -1.5 2.46 -1.7 1.94 -1.7 1.94 -1.5 /
\put{$\bf C(\ )$} at 2.2 -1.9
\endpicture
}\vskip.2in
\centerline{\bf n  = 7, odd type}\vskip.2in
\centerline{
\beginpicture
\setcoordinatesystem units <1in,1in> point   at 0 0
\setplotsymbol({\small.})
\put{$\bf k=3$} at -3 .3
\put7 at -3 0
\put{34} at -3.16 -.35
\plot -3.27 -.25 -3.05 -.25 -3.05 -.45 -3.27 -.45 -3.27 -.25 /
\put{52} at -2.84 -.35
\plot -2.95 -.25 -2.73 -.25 -2.73 -.45 -2.95 -.45 -2.95 -.25 /
\put{322} at -3 -.7
\plot -3.14 -.6 -2.86 -.6 -2.86 -.8 -3.14 -.8 -3.14 -.6 /
\plot -3.08 .1 -2.92 .1 -2.92 -.1 -3.08 -.1 -3.08 .1 /
\plot -3.03 -.1 -3.16 -.25 /
\plot -2.97 -.1 -2.84 -.25 /
\plot -3.16 -.45 -3.04 -.6 /
\plot -2.84 -.45 -2.96 -.6 /
\put{$\bf C(1,3,5,7)$} at -3 -1
\plot -2.5 .4 -2.5 -1.1 /
\put{$\bf k=2$} at -.1 .2
\put{115} at -2.1 -.1
\plot -2.24 0 -1.96 0 -1.96 -.2 -2.24 -.2 -2.24 0 /
\put{1132} at -2.1 -.5
\plot -2.27 -.4 -1.92 -.4 -1.92 -.6 -2.27 -.6 -2.27 -.4 /
\plot -2.1 -.2 -2.1 -.4 /
\put{$\bf C(1.5,7)$} at -2.1 -.8
\put{151} at -1.3 -.1
\plot -1.44 0 -1.16 0 -1.16 -.2 -1.44 -.2 -1.44 0 /
\put{1321} at -1.3 -.5
\plot -1.47 -.4 -1.12 -.4 -1.12 -.6 -1.47 -.6 -1.47 -.4 /
\plot -1.3 -.2 -1.3 -.4 /
\put{$\bf C(1.4,6)$} at -1.3 -.8
\put{511} at -.5 -.1
\plot -.64 0 -.36 0 -.36 -.2 -.64 -.2 -.64 0 /
\put{3211} at -.5 -.5
\plot -.69 -.4 -.32 -.4 -.32 -.6 -.69 -.6 -.69 -.4 /
\plot -.5 -.2 -.5 -.4 /
\put{$\bf C(1,3,5)$} at -.5 -.8
\put{133} at .3 -.1
\plot .16 0 .44 0 .44 -.2 .16 -.2 .16 0 /
\put{1312} at .3 -.5
\plot .11 -.4 .48 -.4 .48 -.6 .11 -.6 .11 -.4 /
\plot .3 -.2 .3 -.4 /
\put{$\bf C(1,4,7)$} at .3 -.8
\put{313} at 1.1 -.1
\plot .96 0 1.24 0 1.24 -.2 .96 -.2 .96 0 /
\put{3112} at 1.1 -.5
\plot .91 -.4 1.28 -.4 1.28 -.6 .91 -.6 .91 -.4 /
\plot 1.1 -.2 1.1 -.4 /
\put{$\bf C(1,3,7)$} at 1.1 -.8
\put{331} at 1.9 -.1
\plot 1.76 0 2.04 0 2.04 -.2 1.76 -.2 1.76 0 /
\put{3121} at 1.9 -.5
\plot 1.71 -.4 2.08 -.4 2.08 -.6 1.71 -.6 1.71 -.4 /
\plot 1.9 -.2 1.9 -.4 /
\put{$\bf C(1,3,6)$} at 1.9 -.8
\put{$\bf k=1$} at -1 -1.1
\put{31111} at -2.4 -1.4
\plot -2.62 -1.3 -2.18 -1.3 -2.18 -1.5 -2.62 -1.5 -2.62 -1.3 /
\put{$\bf C(1,3)$} at -2.4 -1.7
\put{13111} at -1.7 -1.4
\plot -1.92 -1.3 -1.48 -1.3 -1.48 -1.5 -1.92 -1.5 -1.92 -1.3 /
\put{$\bf C(1,4)$} at -1.7 -1.7
\put{11311} at -1 -1.4
\plot -1.22 -1.3 -.78 -1.3 -.78 -1.5 -1.22 -1.5 -1.22 -1.3 /
\put{$\bf C(1,5)$} at -1 -1.7
\put{11131} at -.3 -1.4
\plot -.52 -1.3 -.08 -1.3 -.08 -1.5 -.52 -1.5 -.52 -1.3 /
\put{$\bf C(1,6)$} at -.3 -1.7
\put{11113} at .4 -1.4
\plot .18 -1.3 .62 -1.3 .62 -1.5 .18 -1.5 .18 -1.3 /
\put{$\bf C(1,7)$} at .4 -1.7
\plot .9 -1 .9 -1.8 /
\put{$\bf k=0$} at 1.4 -1.1
\put{1111111} at 1.4 -1.4
\plot 1.1 -1.3 1.7 -1.3 1.7 -1.5 1.1 -1.5 1.1 -1.3 /
\put{$\bf C(1)$} at 1.4 -1.7
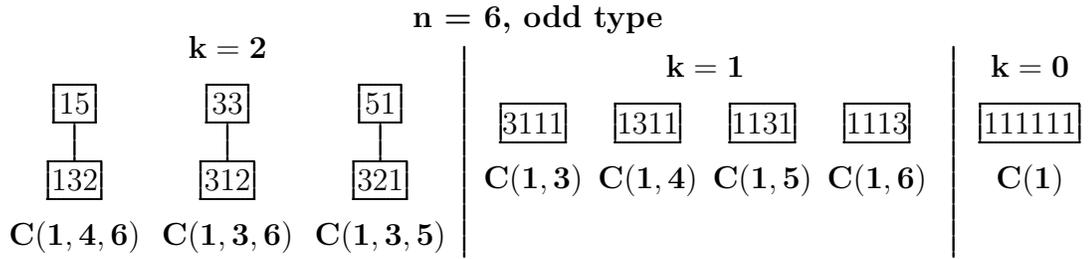
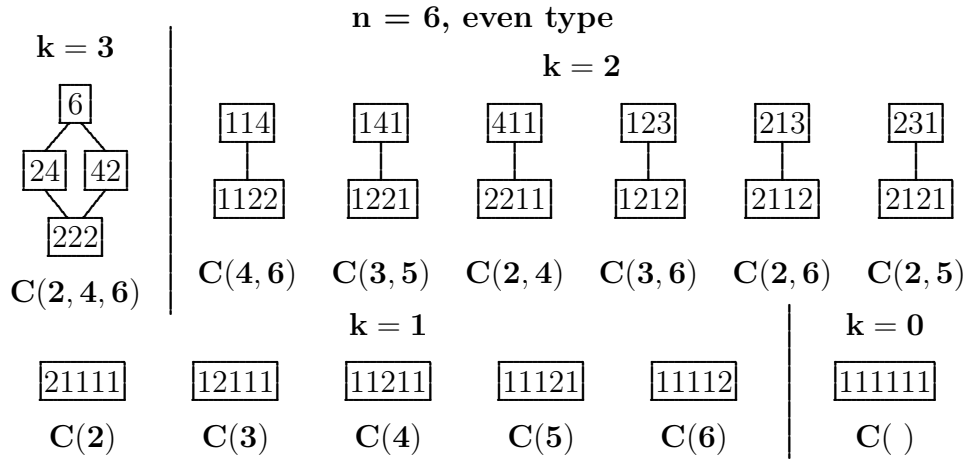
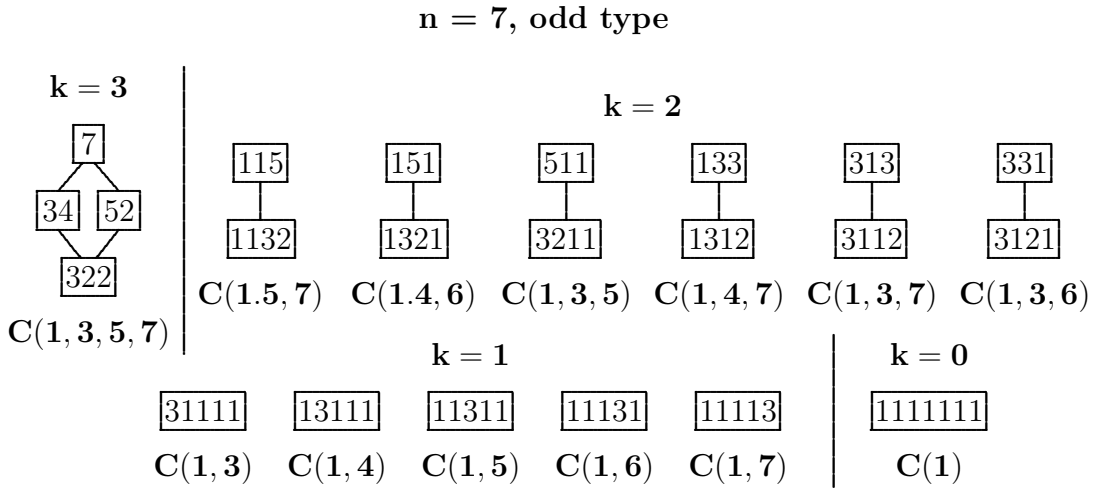
\endpicture
}\vfil\eject

\centerline{\bf n  = 7, even type}\vskip.2in
\centerline{
\beginpicture
\setcoordinatesystem units <1in,1in> point   at 0 0
\setplotsymbol({\small.})
\put{$\bf k=3$} at 0 .3
\put{16} at -1.5 0
\plot -1.61 .1 -1.39 .1 -1.39 -.1 -1.61 -.1 -1.61 .1 /
\put{124} at -1.7 -.35
\plot -1.84 -.25 -1.56 -.25 -1.56 -.45 -1.84 -.45 -1.84 -.25 / 
\put{142} at -1.3 -.35
\plot -1.44 -.25 -1.16 -.25 -1.16 -.45 -1.44 -.45 -1.44 -.25 / 
\put {1222} at -1.5 -.7
\plot -1.67 -.6 -1.32 -.6 -1.32 -.8 -1.67 -.8  -1.67 -.6 /
\plot -1.53 -.1 -1.7 -.25 /
\plot -1.47 -.1 -1.3 -.25 /
\plot -1.7 -.45 -1.54 -.6 /
\plot -1.3 -.45 -1.46 -.6 /
\put{$\bf C(3,5,7)$} at -1.5 -1
\put{25} at -.5 0
\plot -.61 .1 -.39 .1 -.39 -.1 -.61 -.1 -.61 .1 /
\put{214} at -.7 -.35
\plot -.84 -.25 -.56 -.25 -.56 -.45 -.84 -.45 -.84 -.25 / 
\put{232} at -.3 -.35
\plot -.44 -.25 -.16 -.25 -.16 -.45 -.44 -.45 -.44 -.25 / 
\put {2122} at -.5 -.7
\plot -.68 -.6 -.32 -.6 -.32 -.8 -.68 -.8  -.68 -.6 /
\plot -.53 -.1 -.7 -.25 /
\plot -.47 -.1 -.3 -.25 /
\plot -.7 -.45 -.54 -.6 /
\plot -.3 -.45 -.46 -.6 /
\put{$\bf C(2,5,7)$} at -.5 -1
\put{43} at .5 0
\plot .61 .1 .39 .1 .39 -.1 .61 -.1 .61 .1 /
\put{223} at .3 -.35
\plot .44 -.25 .16 -.25 .16 -.45 .44 -.45 .44 -.25 / 
\put{412} at .7 -.35
\plot .84 -.25 .56 -.25 .56 -.45 .84 -.45 .84 -.25 / 
\put {2212} at .5 -.7
\plot .68 -.6 .32 -.6 .32 -.8 .68 -.8  .68 -.6 /
\plot .53 -.1 .7 -.25 /
\plot .47 -.1 .3 -.25 /
\plot .7 -.45 .54 -.6 /
\plot .3 -.45 .46 -.6 /
\put{$\bf C(2,4,7)$} at .5 -1
\put{61} at 1.5 0
\plot 1.61 .1 1.39 .1 1.39 -.1 1.61 -.1 1.61 .1 /
\put{241} at 1.3 -.35
\plot 1.44 -.25 1.16 -.25 1.16 -.45 1.44 -.45 1.44 -.25 / 
\put{421} at 1.7 -.35
\plot 1.84 -.25 1.56 -.25 1.56 -.45 1.84 -.45 1.84 -.25 / 
\put {2221} at 1.5 -.7
\plot 1.67 -.6 1.32 -.6 1.32 -.8 1.67 -.8  1.67 -.6 /
\plot 1.53 -.1 1.7 -.25 /
\plot 1.47 -.1 1.3 -.25 /
\plot 1.7 -.45 1.54 -.6 /
\plot 1.3 -.45 1.46 -.6 /
\put{$\bf C(2,4,6)$} at 1.5 -1
\endpicture
}

\centerline{
\beginpicture
\setcoordinatesystem units <1in,1in> point   at 0 0
\setplotsymbol({\small.})
\put{$\bf k=2$} at 0 -1.3
\put{1114} at -2.7 -1.6
\plot -2.87 -1.5 -2.52 -1.5 -2.52 -1.7 -2.87 -1.7 -2.87 -1.5 /
\put{11122} at -2.7 -2
\plot -2.9 -1.9 -2.49 -1.9 -2.49 -2.1 -2.9 -2.1 -2.9 -1.9 /
\plot -2.7 -1.7 -2.7 -1.9 /
\put{$\bf C(5,7)$} at -2.7 -2.3
\put{1141} at -2.1 -1.6
\plot -2.27 -1.5 -1.92 -1.5 -1.92 -1.7 -2.27 -1.7 -2.27 -1.5 /
\put{11221} at -2.1 -2
\plot -2.3 -1.9 -1.89 -1.9 -1.89 -2.1 -2.3 -2.1 -2.3 -1.9 /
\plot -2.1 -1.7 -2.1 -1.9 /
\put{$\bf C(4,6)$} at -2.1 -2.3
\put{1411} at -1.5 -1.6
\plot -1.67 -1.5 -1.32 -1.5 -1.32 -1.7 -1.67 -1.7 -1.67 -1.5 /
\put{12211} at -1.5 -2
\plot -1.7 -1.9 -1.29 -1.9 -1.29 -2.1 -1.7 -2.1 -1.7 -1.9 /
\plot -1.5 -1.7 -1.5 -1.9 /
\put{$\bf C(3,5)$} at -1.5 -2.3
\put{4111} at -.9 -1.6
\plot -1.07 -1.5 -.72 -1.5 -.72 -1.7 -1.07 -1.7 -1.07 -1.5 /
\put{22111} at -.9 -2
\plot -1.1 -1.9 -.69 -1.9 -.69 -2.1 -1.1 -2.1 -1.1 -1.9 /
\plot -.9 -1.7 -.9 -1.9 /
\put{$\bf C(2,4)$} at -.9 -2.3
\put{1123} at -.3 -1.6
\plot -.47 -1.5 -.12 -1.5 -.12 -1.7 -.47 -1.7 -.47 -1.5 /
\put{11212} at -.3 -2
\plot -.5 -1.9 -.09 -1.9 -.09 -2.1 -.5 -2.1 -.5 -1.9 /
\plot -.3 -1.7 -.3 -1.9 /
\put{$\bf C(4,7)$} at -.3 -2.3
\put{1213} at .3 -1.6
\plot .47 -1.5 .12 -1.5 .12 -1.7 .47 -1.7 .47 -1.5 /
\put{12112} at .3 -2
\plot .52 -1.9 .09 -1.9 .09 -2.1 .52 -2.1 .52 -1.9 /
\plot .3 -1.7 .3 -1.9 /
\put{$\bf C(3,7)$} at .3 -2.3
\put{1231} at .9 -1.6
\plot 1.07 -1.5 .72 -1.5 .72 -1.7 1.07 -1.7 1.07 -1.5 /
\put{12121} at .9 -2
\plot 1.1 -1.9 .69 -1.9 .69 -2.1 1.1 -2.1 1.1 -1.9 /
\plot .9 -1.7 .9 -1.9 /
\put{$\bf C3,6)$} at .9 -2.3
\put{2113} at 1.5 -1.6
\plot 1.67 -1.5 1.32 -1.5 1.32 -1.7 1.67 -1.7 1.67 -1.5 /
\put{21112} at 1.5 -2
\plot 1.72 -1.9 1.29 -1.9 1.29 -2.1 1.72 -2.1 1.72 -1.9 /
\plot 1.5 -1.7 1.5 -1.9 /
\put{$\bf C(2,7)$} at 1.5 -2.3
\put{2131} at 2.1 -1.6
\plot 2.27 -1.5 1.92 -1.5 1.92 -1.7 2.27 -1.7 2.27 -1.5 /
\put{21121} at 2.1 -2
\plot 2.3 -1.9 1.89 -1.9 1.89 -2.1 2.3 -2.1 2.3 -1.9 /
\plot 2.1 -1.7 2.1 -1.9 /
\put{$\bf C(2,6)$} at 2.1 -2.3
\put{2311} at 2.7 -1.6
\plot 2.87 -1.5 2.52 -1.5 2.52 -1.7 2.87 -1.7 2.87 -1.5 /
\put{21211} at 2.7 -2
\plot 2.9 -1.9 2.49 -1.9 2.49 -2.1 2.9 -2.1 2.9 -1.9 /
\plot 2.7 -1.7 2.7 -1.9 /
\put{$\bf C(2,5)$} at 2.7 -2.3
\endpicture
}\vskip.15in
\centerline{
\beginpicture
\setcoordinatesystem units <1in,1in> point   at 0 0
\setplotsymbol({\small.})
\put{$\bf k=1$} at -.75 -2.7
\put{211111} at -2.5 -3
\plot -2.76 -2.9 -2.24 -2.9 -2.24 -3.1 -2.76 -3.1 -2.76 -2.9 /
\put{$\bf C(2)$} at -2.5 -3.3
\put{121111} at -1.8 -3
\plot -2.06 -2.9 -1.54 -2.9 -1.54 -3.1 -2.06 -3.1 -2.06 -2.9 /
\put{$\bf C(3)$} at -1.8 -3.3
\put{112111} at -1.1 -3
\plot -1.36 -2.9 -.84 -2.9 -.84 -3.1 -1.36 -3.1 -1.36 -2.9 /
\put{$\bf C(4)$} at -1.1 -3.3
\put{111211} at -.4 -3
\plot -.66 -2.9 -.14 -2.9 -.14 -3.1 -.66 -3.1 -.66 -2.9 /
\put{$\bf C(5)$} at -.4 -3.3
\put{111121} at .3 -3
\plot .04 -2.9 .56 -2.9 .56 -3.1 .04 -3.1 .04 -2.9 /
\put{$\bf C(6)$} at .3 -3.3
\put{111112} at 1 -3
\plot .74 -2.9 1.26 -2.9 1.26 -3.1 .74 -3.1 .74 -2.9 /
\plot 1.5 -2.6 1.5 -3.4 /
\put{$\bf C(7)$} at 1 -3.3
\put{$\bf k=0$} at 2.1 -2.7
\put{1111111} at 2.1 -3
\plot 1.8 -2.9 2.4 -2.9 2.4 -3.1 1.8 -3.1 1.8 -2.9 /
\put{$\bf C(\ )$} at 2.1 -3.3
\endpicture
}\vskip.15in

Notice that the partition $n=1+1+\dots+1$ belongs to both even and odd types and is considered as two flocks. Notice also that the numbers of partitions of the odd and even type (8 and 13 for $n=6$ and 13 and 21 for $n=7$) agree with the computations in the end of Section 5.3.2.

In the diagrams above, each flock is accompanied by the indication of the container, which is assigned to this flock; we will explain a way to determine this container in the next section. 

\subsubsection{The correspondence containers $\longleftrightarrow$ flocks.}

For every flock, we need to assign a container $C(I)$ with some sparse sequence $I=\{i_1,\dots,i_m\}$. We do this in the following way. 

First of all, to flocks of the odd (even) type, we will assign a sparse sequence, which contains 1 (which does not contain 1). 

If the flock is of even type, then its closing partitions consists of 1's and 2's. More precisely, there are $n-2k$ 1's and $k$ 2's. Let $i_1,\dots,i_k$ be the numbers of parts equal to 2 (we assume that $i_1<\dots<i_k$). Then the container corresponding to this flock is $C(i_1+1,\dots,i_k+k)$. (Notice that the total number of parts is $(n-2k)+k=n-k$, so $m_k+k\le n$.)

If the flock is of odd type, then in its closing partition the first part greater than 1 is 3, and all the subsequent parts are 1's and 2's. Thus, in this partition $n-2k-1$ parts are 1's, $k-1$ parts are 2's, and one part is 3. Let $1_1$ be the number of the part 3, and $i_2,\dots, i_k$ are numbers of parts equal to 2 (again we assume that $m_1<\dots<m_k$). The container corresponding to this flock will be $C(1,1_1+2,1_2+3,\dots,1_k+(k+1))$.(Notice that the total number of parts is $(n-2k-1)+(k-1)+1=n-k$-1, so $m_k+(k+1)\le n$.)

(We should remark that the partition $1+1+\dots+1$ appears as both a flock of the odd type and a flock of the even type. As such, it is assigned to the container $C(1)$  and to the container $C(\ )$.)

It is easy to check that the choice of containers in the cases $n=6$ and $n=7$ shown in the diagrams in Section 5.3 agrees with the rules described here.

\subsection{Characters of the representations in a flock.}\smallskip

The main idea is that for partitions from the same flock, the characters of corresponding irreducible representations look almost the same. We will explain this on the examples considered in Section 5.3.2.  \smallskip

\subsubsection{The case of the even type.} We begin with the example of the even type flock with the head $11=1+4+3+3$ and the tail $11=1+2+2+1+2+1+2\ (n=11, k=4)$ from Section  5.3.3. Let us supplement the dots/dividers diagram corresponding to this flock with the notations $\alpha_1,\alpha_2,\dots, \alpha_{11}$ (corresponding to the dots) and $\beta_1,\beta_2,\dots,\beta_{10}$ (corresponding to the spaces between the dots):

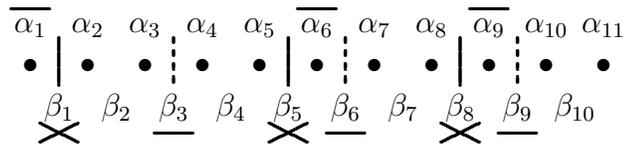
\centerline{
\beginpicture
\setcoordinatesystem units <1in,1in> point   at 0 0
\setplotsymbol(.)
\put{$\alpha_1$} at -1.5 0
\put{$\alpha_2$} at -1.2 0
\put{$\alpha_3$} at -.9 0
\put{$\alpha_4$} at -.6 0
\put{$\alpha_5$} at -.3 0
\put{$\alpha_6$} at 0 0
\put{$\alpha_7$} at .3 0
\put{$\alpha_8$} at .6 0
\put{$\alpha_9$} at .9 0
\put{$\alpha_{10}$} at 1.2 0
\put{$\alpha_{11}$} at 1.5 0
\plot -1.6 .1 -1.4 .1 /
\plot -.1 .1 .1 .1 /
\plot .8 .1 1 .1 /
\put{$\bullet$} at -1.5 -.2
\put{$\bullet$} at -1.2 -.2
\put{$\bullet$} at -.9 -.2
\put{$\bullet$} at -.6 -.2
\put{$\bullet$} at -.3 -.2
\put{$\bullet$} at 0 -.2
\put{$\bullet$} at .3 -.2
\put{$\bullet$} at .6 -.2
\put{$\bullet$} at .9 -.2
\put{$\bullet$} at 1.2 -.2
\put{$\bullet$} at 1.5 -.2
\plot -1.35 -.05 -1.35 -.3 /
\plot -.15 -.05 -.15 -.3 /
\plot .75 -.05 .75 -.3 /
\setdashpattern<2pt,3pt>
\plot -.75 -.05 -.75 -.3 /
\plot .15 -.05 .15 -.3 /
\plot 1.05 -.05 1.05 -.3 /
\setsolid
\put{$\beta_1$} at -1.35 -.42
\put{$\beta_2$} at -1.05 -.42
\put{$\beta_3$} at -.75 -.42
\put{$\beta_4$} at -.45 -.42
\put{$\beta_5$} at -.15 -.42
\put{$\beta_6$} at .15 -.42
\put{$\beta_7$} at .45 -.42
\put{$\beta_8$} at .75 -.42
\put{$\beta_9$} at 1.05 -.42
\put{$\beta_{10}$} at 1.35 -.42
\plot -1.45 -.5 -1.25 -.6 /
\plot -1.45 -.6 -1.25 -.5 /
\plot -.85 -.55 -.65 -.55 /
\plot -.25 -.5 -.05 -.6 /
\plot -.25 -.6 -.05 -.5 /
\plot .05 -.55 .25 -.55 /
\plot .65 -.5 .85 -.6 /
\plot .65 -.6 .85 -.5 /
\plot .95 -.55 1.15 -.55 /
\endpicture
}\vskip.3in

The characters of representations corresponding to the partitions from this flock were computed in Section 4.2. The formula begins with the factor $q^k=q^4$, contain the factors ${\bf e}(x_i\alpha_i)$ for for first $\alpha$ at every odd term in the head partition ($\alpha_1,\alpha_6,$ and $\alpha_7$ in our example, overlined in the diagram and ${\bf e}(y_j\beta_j$ for all $\beta$'s except those corresponding to the solid dividers (marked with crosses in the diagram). Also, in the diagram we underlined the $\beta$'s, which correspond to the dotted dividers; they may disappear in the characters of representations, corresponding to non-head partitions in the flocks. It is convenient to assume that the underlined $y$'s are allowed to be zero, while the terms not underlined must be different from zero. So, here is the character:$$q^4{\bf e}(x_1\alpha_1){\bf e}(x_6\alpha_6){\bf e}(x_9\alpha_9){\bf e}(y_2\beta_2)\underline{{\bf e}(y_3\beta_3)}{\bf e}(y_4\beta_4)\underline{{\bf e}(y_6\beta_6)}{\bf e}(y_7\beta_7)\underline{{\bf e}(y_9\beta_9)}{\bf e}(y_{10}\beta_{10})$$if $\alpha_2=\alpha_3=\alpha_4=\alpha_7=\alpha_{10}=0$ and $\alpha_6\colon\alpha_8=y_7\colon y_6,\alpha_9\colon\alpha_{11}=y_{10}\colon y_9$, and is zero otherwise. It is worth noticing that we underline {\it every second} term with $\beta$.

The transition from the example to the general case is straightforward. We will only add the formulas for the3 numbers of the terms with $\alpha$'s and underlined and not underlined terms with $\beta$'s.

The number of overlined $\alpha$' s (same as the number of factors with $\alpha$'s in the formula) is the number of odd terms in the partitions of the flock: it is $n-2k$. 

The number of crossed $\beta$'s is the number of solid dividers, which is one less that the number of terms in the head partition; so, it is also $n-2k$. Thus the number of $\beta$'s in the formula for the character is $(n-1)-(n-2k)=2k-1$.

The number of underlined $\beta$-factors is the number of dotted dividers is $k-1$ (see Section 5.3.3). Hence, the number of $\beta$-factors not underlined is $(2k-1)-(k-1)=k$.

\subsubsection{The case of the odd type.} Again, we begin with the example of the odd type flock from Section  5.3.3. The head and the tail of this flock are $12=3+1+5+3$ and $12=3+1+1+2+2+1+2$ (thus, $n=12,k=4$). We consider the dots/dividers diagram similar to that in Section 5.4.1.\vskip.2in

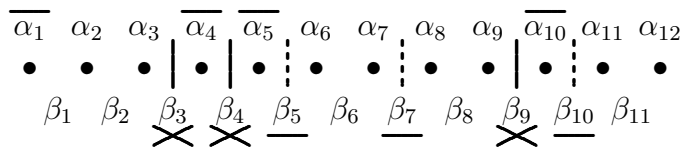
\centerline{
\beginpicture
\setcoordinatesystem units <1in,1in> point   at 0 0
\setplotsymbol(.)
\put{$\alpha_1$} at -1.65 0
\put{$\alpha_2$} at -1.35 0
\put{$\alpha_3$} at -1.05 0
\put{$\alpha_4$} at -.75 0
\put{$\alpha_5$} at -.45 0
\put{$\alpha_6$} at -.15 0
\put{$\alpha_7$} at .15 0
\put{$\alpha_8$} at .45 0
\put{$\alpha_9$} at .75 0
\put{$\alpha_{10}$} at 1.05 0
\put{$\alpha_{11}$} at 1.35 0
\put{$\alpha_{12}$} at 1.65 0
\plot -1.75 .1 -1.55 .1 /
\plot -.85 .1 -.65 .1 /
\plot -.55 .1 -.35 .1 /
\plot .95 .1 1.15 .1 /
\put{$\bullet$} at -1.65 -.2
\put{$\bullet$} at -1.35 -.2
\put{$\bullet$} at -1.05 -.2
\put{$\bullet$} at -.75 -.2
\put{$\bullet$} at -.45 -.2
\put{$\bullet$} at -.15 -.2
\put{$\bullet$} at .15 -.2
\put{$\bullet$} at .45 -.2
\put{$\bullet$} at .75 -.2
\put{$\bullet$} at 1.05 -.2
\put{$\bullet$} at 1.35 -.2
\put{$\bullet$} at 1.65 -.2
\plot -.9 -.05 -.9 -.3 /
\plot -.6 -.05 -.6 -.3 /
\plot .9 -.05 .9 -.3 /
\setdashpattern<2pt,3pt>
\plot -.3 -.05 -.3 -.3 /
\plot .3 -.05 .3 -.3 /
\plot 1.2 -.05 1.2 -.3 /
\setsolid
\put{$\beta_1$} at -1.5 -.42
\put{$\beta_2$} at -1.2 -.42
\put{$\beta_3$} at -.9 -.42
\put{$\beta_4$} at -.6 -.42
\put{$\beta_5$} at -.3 -.42
\put{$\beta_6$} at 0 -.42
\put{$\beta_7$} at .3 -.42
\put{$\beta_8$} at .6 -.42
\put{$\beta_9$} at .9 -.42
\put{$\beta_{10}$} at 1.2 -.42
\put{$\beta_{11}$} at 1.5 -.42
\plot -1 -.5 -.8 -.6 /
\plot -1 -.6 -.8 -.5 /
\plot -.7 -.5 -.5 -.6 /
\plot -.7 -.6 -.5 -.5 /
\plot -.4 -.55 -.2 -.55 /
\plot .2 -.55 .4 -.55 /
\plot .8 -.5 1 -.6 /
\plot .8 -.6 1 -.5 /
\plot 1.1 -.55 1.3 -.55 /
\endpicture
}\vskip.3in

The formula for the character is obtain according to he same rules as in the even case: $$\begin{array} {l} q^4{\bf e}(x_1\alpha_1){\bf e}(x_4\alpha_4){\bf e}(x_5\alpha_5){\bf e}(x_{10}\alpha_{10})\cdot\\ \hskip1.2in{\bf e}(y_1\beta_1){\bf e}(y_2\beta_2)\underline{{\bf e}(y_5\beta_5)}{\bf e}(y_6\beta_6)\underline{{\bf e}(y_7\beta_7)}{\bf e}(y_8\beta_8)\underline{{\bf e}(y_{10}\beta_{10})}{\bf e}(y_{11}\beta_{11}),\end{array} \eqno(11)$$if $\alpha_2=\alpha_6=\alpha_8=\alpha_{11}=0$ and $\alpha_3\colon\alpha_1=y_1\colon y_2, \alpha_7\colon\alpha_5=y_5\colon y_6, \alpha_9\colon\alpha_7=y_7\colon y_8, \alpha_{12}\colon\alpha_{10}=y_{10}\colon y_{11}$, and is zero otherwise.

The number of overlined $\alpha$'s is the same as in the even case: $n-2k$ (it is the number of odd terms in the partition).

The number of crossed $\beta$'s is again the number of solid dividers, thus it it one less than the number of terms in the partition. Bur now it is $n-2k$, so the number of crossed $\beta$'s is $n-2k-1$. Thus the number of $\beta$'s in the formula for the characters is $(n-1)-(n-2k-1)=2k$ .

The number of underlined $\beta$-factors is again the number of dotted dividers is $k-1$. Hence, the number of $\beta$-factors not underlined is $2k-(k-1)=k+1$.

\subsection{ The end of the construction of a model.} It remains to specify the one-dimensional representations (the characters) of the stabilizer of the classes within the container corresponding to a flock. Again, we will give the details of the construction for the examples of Section 5.3.3 and then discuss the general case.\smallskip

\subsubsection{The case of the even type; the example.} Again, we begin with the example of the flock with the head $11=1+4+3+3$ and the tail $11=1+2+2+1+2+1+2$ from Section 5.3.3. According to Section 5.3.5, the container corresponding to this section is $C(I)$, where $I=\{3,5,8,11\}$. In the notations of Section 3.3.2, $I^-=\{2,4,6,7,9,10\}$ and $I^+=\{1,3,5,8,11\}$. The classes in this container are characterized by\smallskip

4 $a$-invariants, $a_3,a_5,a_8,a_{11}\in{\mathcal F}_q-0$ and

4 $b$-invariants, $b_	1,b_3a_5+b_4a_3, b_6, b_9\in {\mathcal F}_q$. \smallskip

The characters of representations induced by 1-dimensional representations of the common stabilizer 
$$\Stab(i)=\{g(\alpha_1,0,\alpha_3,0,\alpha_5,0,0,\alpha_8,0,0,\alpha_{11};\beta_1,\dots,\beta_{10})\}$$
\noindent of classes in $C(I)$ are $$(\overline\alpha,\overline\beta)_\mapsto q^6{\bf e}(A_1\alpha_1)(A_3\alpha_3)(A_5\alpha_5)(A_8\alpha_8)(A_{11}\alpha_{11}){\bf e}(B_1\beta_1)\dots{\bf e}(B_{10}\beta_{10}),\eqno(10)$$if for $i\in I^-,$ $\alpha_i=0$ and $B_{i-1}\alpha_{i-1}=B_i\alpha_{i+1}$ and 0 otherwise.

The characters of representations corresponding to our flock are $$q^4{\bf e}(x_1\alpha_1){\bf e}(x_6\alpha_6){\bf e}(x_9\alpha_9){\bf e}(y_2\beta_2)\underline{{\bf e}(y_3\beta_3)}{\bf e}(y_4\beta_4)\underline{{\bf e}(y_6\beta_6)}{\bf e}(y_7\beta_7)\underline{{\bf e}(y_9\beta_9)}{\bf e}(y_{10}\beta_{10}),\eqno(11)$$where $x$'s and $y$'s are element of ${\mathbb F}_q$, different from 0 for $y$'s not underlined. We need to make choose values for $A$'s and $B$'s in terms of $a$- and $b$-invariants to accommodate the formulas (10) and (11). For values of $B_2,B_4,B_7,B_{10}$ we take (in any order) the $a$-invariants $a_3,a_5,a_8,a_{11}$. For $B_3,B_6, B_9$, we take (in any order) three of the four $b$-invariants, say, $b_3a_5+b_4a_3,b_6,b_9$. the remaining $b$-invariant, $b_1$, we take for the value of $A_1$ (taking into consideration the fact that $\alpha_1$ appears in both formulas (10) and (11)). For the remaining $A_i$'s and $B_j$'s in (10) we assume zero values.

Still two discrepancies in the formulas (10) and (11) remain. First, the  first factors in (10) and (11) are different: $q^6$ and $q^4$. Second, (11) contains terms with $\alpha_6$ and $\alpha_9$, while (10) does not. Both can be eliminated by the following move. We remove the conditions $\alpha_6=\alpha_9=0$ and replace $q^6=q^4\cdot q^2$ by $q^4\cdot\left(\sum\limits_{A_6\in{\mathbb F}_q}A_6\alpha_6\right)\cdot\left(\sum\limits_{A_9\in{\mathbb F}_q}A_9\alpha_9\right).$ Then the character (10) (of capacity $q^6$) becomes the sum of characters (11) of $q^4$-dimensional representations labelled by $A_6,A_9\in{\mathbb F}_q$.    

\subsubsection{The case of the odd type; the example.} Now we consider the flock with the  with the head $12=3+1+5+3$ and the tail $12=3+1+1+2+2+1+2$ (see Section 5.3.3). According to Section 5.3.5, the container corresponding to this section is $C(I)$, where $I=\{1,3,7,9,12\}$. In the notations of Section 3.3.2, $I^-=\{2,4,6,8,10,11\}$ and $I^+=\{1,3,5,7,9,12\}$. The classes in this container are characterized by\smallskip

5 $a$-invariants, $a_1,a_3,a_7,a_9,a_{12}\in{\mathcal F}_q-0$ and

5 $b$-invariants, $b_	1a_1+b_2a_3,b_4,b_5,b_7a_7+b_8a_9,b_{10}\in {\mathcal F}_q$. \smallskip

The formulas similar to (10) and (11) are $$q^6{\bf e}(A_1\alpha_1){\bf e}(A_3\alpha_3){\bf e}(A_5\alpha_5){\bf e}(A_7\alpha_7){\bf e}(A_9\alpha_9){\bf e}(A_{12}\alpha_{12}){\bf e}(B_1\beta_1)\dots{\bf e}(B_{11}\beta_{11})\eqno(12) $$and $$\begin{array} {l} q^4{\bf e}(x_1\alpha_1){\bf e}(x_4\alpha_4){\bf e}(x_5\alpha_5){\bf e}(x_{10}\alpha_{10}) \cdot \\ \hskip1in {\bf e}(y_1\beta_1){\bf e}(y_2\beta_2)\underline{{\bf e}(y_5\beta_5)}{\bf e}(y_6\beta_6)\underline{{\bf e}(y_7\beta_7)}{\bf e}(y_8\beta_8)\underline{{\bf e}(y_{10}\beta_{10})}{\bf e}(y_{11}\beta_{11})\end{array}\eqno(13)$$The coefficients $y$ in the non-underlined  factors in (13) must be non-zero, so we used for them the $a$-invariants: $a_1,a_3,a_7,a_9,a_{12}$ for $B_1,B_2,B_6,B_8,B_{11}$. For $B_5, B_7, B_{10}$, we use the $b$-invariants, say, $b_5,b_7a_7+b_8a_9,b_{10}$, and the remaining $b$-invariants we assign to the coefficients $A_1,A_5$ in the terms with $\alpha_1$ and $\alpha_5$, which appear in the both formulas (12) and (13). All the other coefficients $A$ and $B$ in (12) are assumed to be zero.

The last step repeats the last step in the previous case: we remove the conditions $\alpha_4=\alpha_{10}=0$ and replace $q^6$ in (12) by $q^4\cdot\left(\sum\limits_{A_4\in{\mathbb F}_q}A_4\alpha_4\right)\cdot\left(\sum\limits_{A_{10}\in{\mathbb F}_q}A_{10}\alpha_{10}\right).$
 
\subsubsection{General case} Let us consider ion the general case the construction of the character of the stabilizer which was described above for examples. 

First we take the non-underlined factors ${\bf e}(y\beta)$ in  the formulas like (11) and (13) and assign for their coeficients the $a$-invariants. But the numbers of non-underlined terms and $a$-invariants agree: they both are $k$ for the even type and  $k+1$ for the odd type. Thus, this step works in the general case. 

Next we use $b$-invariants for the underlined factors ${\bf e}(y\beta)$. The number of this factors is $|I|-1$ in the even case and $|I|-2$ in the odd case, while the total number of $b$-invariants is always $|I^+|-1$. But $I^+\supseteq I$, so $|I^+|\geq|I|$, hence there are sufficient quantity of $b$-invariants. Moreover, $|I^+|-|I|$ $b$-invariants in the even case and  $|I^+|-|I|+1$ $b$-invariants in the odd case and remain for the next step.

Next step consists in assigning the remaining $b$-invariants to the coefficients to $\alpha_i$, which appear both characters we consider. We need to check that the numbers of these $b$-invariants and these $\alpha_i$ are the same. The cases of the even and odd types are slightly different. In the {\bf even type case} the number of remaining $b$-invariants is $|I^+-I|$. The elements of $I^+-I$ correspond to those 1's ion the tail partition, which do not follow 2; the same 1's correspond to the common $\alpha's$ in the two characters. The only difference in the {\bf odd type case} is that 1 belongs to $I$, and hence not to $I^+-I$, while $\alpha_1$ appears in the both characters. This leads to the additional 1 in the expression $|I^+-I|+1$ for the number of common $\alpha_i$'s in the two characters.

And the last step: the difference between $|I^-|$ and $k$ in the two characters should be compensated by the unused part of the number of $\alpha_i$ in the character of the representations in the flock $|I^+-I|=|I^+|-k$. In is easy:$$(n-2k)-(|I^+|-k)=n-2k-|I^+|+k=n-2k-(n-|I^-|)+k=|I^-|-k.$$

\subsubsection{The possible non-commutativity of the stabilizer does not matter.} If the stabilizer $\Stab(I)$ is not commutative, then some $\beta$-factors in the formula for the character of the stabilizer (like (10) or (12)) may be missing. But this does not affect our construction. The reason is that this missing $\beta$'s appear, when there are gaps in $I$ of the length 3 or more (even 4 or more, if this gap {\it inside} I, not in one of the ends). In this case, there are at least two 1's in a row in the tail partition of the flock, they have to be separated by a solid divider (there are at most one 1 between any two successive dividers), and the corresponding $\beta$ is missing also in the character of the representation in the flock (formula like (11) or (13)).

\end{document}